%% file: ms.tex
\documentclass[opre,nonblindrev]{oo}

\OneAndAHalfSpacedXI 

\usepackage{endnotes}
\let\footnote=\endnote
\let\enotesize=\normalsize
\def\notesname{Endnotes}%
\def\makeenmark{\hbox to1.275em{\theenmark.\enskip\hss}}
\def\enoteformat{\rightskip0pt\leftskip0pt\parindent=1.275em
  \leavevmode\llap{\makeenmark}}

\input{Definitions}

\usepackage{natbib}
 \bibpunct[, ]{(}{)}{,}{a}{}{,}%
 \def\bibfont{\small}%
 \def\bibsep{\smallskipamount}%
\TheoremsNumberedThrough     
\ECRepeatTheorems

\EquationsNumberedThrough    

\begin{document}


\RUNAUTHOR{Guo, Bodur, and Taylor}

\RUNTITLE{Copositive Duality for Discrete Markets and Games}

\TITLE{Copositive Duality for Discrete Markets and Games}

\ARTICLEAUTHORS{%
\AUTHOR{Cheng Guo, Merve Bodur}
\AFF{Department of Mechanical and Industrial Engineering, University of Toronto, Toronto, ON M5S 3G8, \EMAIL{cguo@mie.utoronto.ca, bodur@mie.utoronto.ca}} 
\AUTHOR{Joshua A. Taylor}
\AFF{The Edward S. Rogers Sr. Department of Electrical and Computer Engineering, University of Toronto, ON M5S 3G4, \EMAIL{josh.taylor@utoronto.ca}}
} 

\ABSTRACT{Optimization problems with discrete decisions are nonconvex and thus lack strong duality, which limits the usefulness of tools such as shadow prices and the KKT conditions. It was shown in \cite{burer2009copositive} that mixed-binary quadratic programs can be written as completely positive programs, which are convex. Completely positive reformulations of discrete optimization problems therefore have strong duality if a constraint qualification is satisfied. We apply this perspective in two ways. First, we write unit commitment in power systems as a completely positive program, and use the dual copositive program to design a new pricing mechanism. Second, we reformulate integer programming games in terms of completely positive programming, and use the KKT conditions to solve for pure strategy Nash equilibria. To facilitate implementation, we also design a cutting plane algorithm for solving copositive programs exactly. 
}



\maketitle

%


\input{content2}

\bibliographystyle{informs2014} 
\begingroup
    \setlength{\bibsep}{-0.5pt}
    \linespread{1}\selectfont
\bibliography{MainBib_josh,JATBib_josh,reference}
\endgroup

\ECSwitch


\ECHead{e-Companion}

\input{ecompanion}

%
%
%






\end{document}

%% file: Definitions.tex
\usepackage{amsmath}
\usepackage{amssymb}
\usepackage{url}
\usepackage[colorlinks = true,
            linkcolor = blue,
            urlcolor  = blue,
            citecolor = blue,
            anchorcolor = blue]{hyperref}
\usepackage{enumitem}
\usepackage[utf8]{inputenc}
\usepackage[english]{babel}
\usepackage[table,xcdraw]{xcolor}
\usepackage{booktabs}
\usepackage[titletoc]{appendix}
\usepackage{subcaption}
\usepackage{bm}
\usepackage{comment}
\usepackage{ulem}
\usepackage{tikz}
\usetikzlibrary{shapes,arrows.meta, arrows,quotes,positioning,tikzmark,fit,backgrounds}
\tikzstyle{block} = [rectangle, draw, text width=10em, text centered, minimum height=2em]
\tikzstyle{box} = [rectangle, draw,text width=5em,text centered]
\tikzstyle{line} = [draw, -latex']
\tikzstyle{block2} = [rectangle, draw, fill=gray!20, text width = 17em, text centered, rounded corners, minimum height=1em]
\tikzstyle{box2} = [rectangle,text width=5em,text centered]

\newcommand{\bsubeq}{\begin{subequations}}
\newcommand{\esubeq}{\end{subequations}}
\newcommand{\BI}{\begin{itemize}}
\newcommand{\EI}{\end{itemize}}
\newcommand{\I}{\item}

\newcommand{\cred}{\color{black}}

\renewcommand{\mc}{\mathcal}
\newcommand{\mb}{\mathbb}
\newcommand{\Pith}{\mc{P}_i}

\newcommand{\diag}{\text{diag}}
\newcommand{\gen}{\text{Gen}}
\newcommand{\MP}{\text{(MP)~}}

\newcommand{\SP}{\mc{SP}}
\newcommand{\bq}{\text{MBQ}}
\newcommand{\bqp}{\text{MBQP}}
\newcommand{\cp}{\text{CP}}
\newcommand{\cpp}{\text{CPP}}
\newcommand{\cop}{\text{COP}}

\newcommand{\uc}{\mathcal{UC}}

\newcommand{\Qi}{Q^{(i)}}
\newcommand{\Bi}{\mc{B}_i}
\renewcommand{\c}{\mathbf{c}^{\top}}
\newcommand{\ci}{\mathbf{c}^{(i)\top}}
\renewcommand{\a}{\mathbf{a}}
\newcommand{\alambda}{\mathbf{h}}
\newcommand{\alambdatop}{\mathbf{h}^\top}
\newcommand{\aphi}{\a}
\newcommand{\aphitop}{\a^{\top}}

\newcommand{\azero}{\a_0}
\newcommand{\azerotop}{\a_0^\top}
\newcommand{\Aone}{A_{1}}
\newcommand{\Aonetop}{A_1^\top}
\newcommand{\Ag}{A_g}
\newcommand{\Agtop}{A_g^\top}
\newcommand{\AG}{A_{|\mathcal{G}|}}
\newcommand{\AGtop}{A_{|\mathcal{G}|}^\top}

\newcommand{\aj}{\mathbf{a}_{j}}
\newcommand{\ajtop}{\mathbf{a}_{j}^\top}
\newcommand{\aij}{\mathbf{a}_j^{(i)}}
\newcommand{\aijtop}{\mathbf{a}_j^{(i)\top}}

\newcommand{\bij}{b_j^{(i)}}
\newcommand{\ek}{\mathbf{e}_k}
\newcommand{\ektop}{\mathbf{e}_k^{\top}}
\newcommand{\el}{\mathbf{e}_l}
\newcommand{\eltop}{\mathbf{e}_l^{\top}}
\newcommand{\y}{\mathbf{y}}

\newcommand{\alphai}{{\pmb\alpha}^{(i)}}
\newcommand{\alphaitop}{{\pmb\alpha}^{(i)\top}}

\newcommand{\p}{\mathbf{p}}
\newcommand{\x}{\mathbf{x}}

\newcommand{\xstarik}{\mathbf{x}_{ik}^*}
\newcommand{\xg}{\mathbf{x}_g}
\newcommand{\xstarg}{\mathbf{x}_g^{*}}
\newcommand{\xstar}{\mathbf{x}^*}

\newcommand{\xith}{\mathbf{x}_i}
\newcommand{\xitop}{\mathbf{x}_i^{\top}}
\newcommand{\xnitop}{\mathbf{x}_{-i}^{\top}}
\newcommand{\xni}{\mathbf{x}_{-i}}
\newcommand{\xstari}{\mathbf{x}_i^*}
\newcommand{\xstaritop}{\mathbf{x}_i^{*\top}}
\newcommand{\xstarni}{\mathbf{x}_{-i}^*}

\newcommand{\Xith}{X_i}
\newcommand{\Xp}{X^{pp}_{gt, gt}}
\newcommand{\Xpstar}{X^{pp*}_{gt, gt}}
\newcommand{\Xpp}{X^{pp}_{gt, g't}}
\newcommand{\Xppstar}{X^{pp*}_{gt, g't}}
\newcommand{\Xvw}{X^{vw}_{g_1 t, g_2 t}}
\newcommand{\Xvwstar}{X^{vw*}_{g_1 t, g_2 t}}
\newcommand{\mcXi}{\mc{X}_i}
\newcommand{\mcXip}{\mc{X}^{'}_{i}}
\newcommand{\mcXni}{\mc{X}_{-i}}
\newcommand{\Xstari}{X_i^*}
\newcommand{\opt}{\text{opt}}

\newcommand{\feas}{\text{Feas}}
\newcommand{\tr}{\text{Tr}}
\newcommand{\xhat}{\hat{\mathbf{x}}}
\newcommand{\xbar}{\bar{\x}}
\newcommand{\Xhat}{\hat{X}}
\newcommand{\Xbar}{\bar{X}}

\def\st{{\rm s.t.}}

\usepackage[ruled,vlined]{algorithm2e}

\newcommand{\conv}{{\rm conv}} 

\usepackage{longtable}
\newcommand{\sm}{\setminus}

\usepackage{mathtools}

\usepackage{multirow}
\usepackage{graphicx}
\usepackage{caption}
\usepackage{breqn}
\allowdisplaybreaks

%% file: content2.tex
\section{Introduction}
Discreteness is a challenge in many contexts. Examples include optimization, markets with discrete decisions, which can lack efficient equilibria, and games. A basic difficulty is the lack convexity, which precludes the use of tools like strong duality and the Karush-Kuhn-Tucker (KKT) conditions. Many discrete problems can be written as mixed-integer programs (MIPs). \cite{burer2009copositive} has shown that mixed-binary quadratic programs (MBQPs), a generalization of MIPs, can be written as completely positive programs (CPPs). CPPs and their dual copositive programs (COPs) are convex but NP-hard. While this does not point to a better way of solving MBQPs, it does provide a new notion of duality. In this paper, we explore the use of copositive duality for discrete problems, focusing on two applications: pricing in nonconvex electricity markets, and characterizing the equilibria of discrete games.


An immediate challenge in using copositive duality is the relatively small number of options for solving COPs exactly. Our first contribution is the design of a novel cutting plane algorithm that exactly solves COPs when it terminates. The algorithm consists of a sequence of MIPs, which can thus be implemented using standard industrial solvers. It also accommodates discrete variables, and can thus solve COPs with integer constraints. This feature is the first in the literature and useful for obtaining the equilibria of integer programming games (IPGs) in Section \ref{ch: app_game}. 

Our first application of copositive duality is unit commitment (UC) in power systems, in which the decision to turn a generator on or off is binary. It is commonly formulated as an MIP, for which efficient prices are difficult to construct; see,  e.g., \cite{liberopoulos2016critical}. We rewrite the MIP as a CPP, and use the dual COP to design a pricing mechanism based on shadow prices. The mechanism is budget balanced and, under certain conditions, individually rational. It is also straightforward to incorporate additional features, such as revenue adequacy of the generators, by adding constraints directly to the dual COP.

In our second application, we use copositive duality to characterize the Nash equilibria (NE) of IPGs with both binary and continuous decisions. There are well-known conditions guaranteeing the existence and uniqueness of NE in convex (concave) games with convex strategy sets~\citep{debreu1952social,fan1952fixed, glicksberg1952further, rosen1965existence}. Also, due to the convexity of strategy sets, the KKT conditions can be used to compute the NE. We use copositive duality to extend some these results to IPGs, and use the KKT conditions to compute their equilibria.

The reminder of the paper is organized as follows. In Section \ref{ch: lit_review} we review the literature on copositive programming, pricing in nonconvex markets, and discrete games. In Section \ref{ch: cop_duality} we 
provide the necessary background on CPPs and COPs. In Section \ref{ch: cut_plane} we present the new cutting plane algorithm for solving COPs. 
In Section \ref{ch: app_pricing} we design a COP-based pricing scheme for UC, and in Section \ref{ch: app_game} we derive new results on the equilibria of IPGs. In Section \ref{ch: results} we present the results of our computational experiments. We conclude the paper in Section \ref{sec:conc}.

All proofs are in Section \ref{ec: proofs} of the e-companion.

\section{Literature Review}\label{ch: lit_review}
In this section, we review the relevant literature. Section \ref{ch: review_cop} reviews the relationship between CPPs, COPs, and integer programs, and solution methods for COPs. Section \ref{ch: review_pricing} surveys the pricing schemes for nonconvex markets. Section \ref{ch: review_game} reviews on the literature on IPGs.
\subsection{Copositive Programming}\label{ch: review_cop}
Copositive programming has been shown to generalize a number of NP-hard problems, such as quadratic optimization \citep{bomze2002solving}, two-stage adjustable robust optimization \citep{xu2018copositive, hanasusanto2018conic}, and MBQPs \citep{burer2009copositive}. In particular, \cite{burer2009copositive} shows that an MBQP can be written as a CPP. This reformulation is the basis of our work. 

At present, no industrial software can directly solve COPs. \cite{parrilo2000structured} constructs a hierarchy of semidefinite programs (SDPs), which is widely used to approximate COPs. For example, it is used by \cite{de2002approximation} 
to find the stability number of a graph, and by 
\cite{hanasusanto2018conic} 
for two-stage distributionally robust linear programs.  
An exact algorithm for COPs based on simplicial partitions is proposed by \cite{bundfuss2009adaptive}. \cite{bomze2008new} and \cite{bomze2010copositivity} use cutting planes to strengthen the SDP relaxation for COPs, but both are problem specific and can only be used for quadratic programming and clique number problems, respectively. In this paper, we develop a purely MIP-based cutting plane algorithm that exactly solves COPs when it terminates. Also, to the best of our knowledge, this is the first algorithm that can solve COPs with discrete variables.
\subsection{Pricing in Discrete Markets}\label{ch: review_pricing}
Many markets have nonconvexities, e.g., due to binary decisions and indivisible goods. It is difficult to design efficient pricing mechanisms in such markets because of the lack of strong duality. This has now been a subject of research for decades. We refer the reader to \cite{liberopoulos2016critical} for a more thorough review of pricing in nonconvex markets.

A basic source of nonconvexity in electricity markets is the binary startup and shutdown decisions of generators. These decisions are optimized via UC, which is commonly formulated as an MIP \citep{carrion2006computationally}. The basic idea of most current approaches is to construct approximate shadow prices for the MIP. \cite{o2005efficient} eliminate the nonconvexity by fixing the binary decisions at their optimal values and pricing the resulting linear program (LP). We call this scheme restricted pricing (RP). RP and its variants are used by some independent system operators (ISOs) in the US, such as the Pennsylvania-New Jersey-Maryland Interconnect. However, RP is often too low to cover the costs of generators, in which uplift or make-whole payments are needed to ensure generator profitability. The  convex hull pricing (CHP) of \cite{hogan2003minimum} and \cite{gribik2007market} uses the Lagrangian multipliers of the demand constraints as prices, and has been shown shown to minimize uplift. A modified version of CHP called extended locational marginal pricing is used by the Midcontinent ISO. \cite{ruiz2012pricing} propose a primal-dual approach for pricing, which combines the UC problem and the dual of its linear relaxation, as well as revenue adequacy constraints that ensure nonnegative profit for each generator. We note that such revenue adequacy constraints eliminate uplift, but also modify the original UC problem and may therefore result in suboptimal decisions. 

We use the CPP reformulation of MBQP to design a new pricing scheme, which we refer to as copositive duality pricing (CDP). CDP is budget balanced, supports the optimal UC, and is flexible in that it provides direct access to the dual COP. We make use of this flexibility by adding a revenue adequacy constraint. We refer to the augmented pricing scheme as revenue-adequate CDP (RCDP). 

It is desirable for a pricing scheme to ensure individual rationality, i.e., that individual generators have no incentive to deviate from the optimal UC solution. \cite{o2005efficient} prove that RP is individually rational. CHP satisfies individual rationality in some special cases, but in general it does not support each individual generator's profit-maximizing solution \citep{gribik2007market}. The primal-dual approach does not guarantee individual rationality either. CDP does not in general lead to individual rationality, but we provide a simple sufficient condition under which it does.

Table \ref{tab:price_features} compares the properties of several pricing schemes, where ``$\circ$" indicates the existence of a sufficient condition in the literature for ensuring the property. 
\begin{table}[h]
\centering
\caption{Features of Pricing Schemes}
\label{tab:price_features}
\begin{tabular}{lccc}
\hline
\multicolumn{1}{c}{} & Uplift & Optimal & Individual  \\
Scheme               & free  & UC      & rationality \\ \hline
RP                  & $\times$    & $\checkmark$     & $\checkmark$         \\
CHP                  & $\times$    & $\times$      & $\circ$          \\
Primal-dual          & $\checkmark$    & $\times$      & $\times$          \\
CDP                  & $\times$    & $\checkmark$     & $\circ$     \\
RCDP         & $\checkmark$     & $\times$     & $\times$         \\ \hline
\end{tabular}
\end{table}

We also mention that there is a related literature stream focusing on markets with indivisible goods. Recent papers include \cite{danilov2001discrete} and \cite{baldwin2019understanding}, which use discrete convexity to prove the existence of equilibria.

\subsection{Integer Programming Games}\label{ch: review_game}
IPGs are a class of games where each player's action set contains integer decisions. Some examples of IPGs include bimatrix games, Nash-Cournot energy production games \citep{gabriel2013solving}, and Cournot games with indivisible goods \citep{kostreva1993combinatorial}. Most papers focus on algorithms for computing NE, and are limited to IPGs with only integer decisions \citep{koppe2011rational, sagratella2016computing, carvalho2017computing}. Theoretical results on IPGs include \cite{mallick2011existence} and \cite{sagratella2016computing}, which respectively provide NE existence conditions for two-person discrete games and 2-groups partitionable discrete games. In this work we provide conditions for the existence and uniqueness of pure NE (PNE) in IPGs with both continuous and discrete variables, as well as KKT conditions that can be used to compute the PNE. 

Similar to our approach, in the literature there are several applications of convex optimization to IPGs. \cite{parrilo2006polynomial} reformulates zero-sum polynomial games as a single SDP. \cite{ahmadi2020semidefinite} use SDP to find the additive $\epsilon$-approximate NE of bimatrix games. Closely related to our work, \cite{sayin2019optimal} study a Stackelberg game with binary decisions in the context of optimal hierarchical signaling. They reformulate each player's optimization problem as a CPP, and solve the CPP and its dual COP approximately with SDP. Our work differs in that we analyze the equilibria of general IPGs, and for certain classes we can obtain the PNE by solving a single COP exactly.
\section{Background}\label{ch: cop_duality}
In this section, we define CPP and COP and state some basic results. 

Throughout the paper, we use bold letters for vectors. $\tr(\cdot)$ denotes the trace of a matrix and $(\cdot)^\top$ denotes the transpose of a vector or matrix. $\ek\in\mb{R}^n$ is the $k^{\text{th}}$ unit vector.
\subsection{Preliminaries}\label{ch: preliminary}
Let $\mc{S}_n$ be the set of $n$-dimensional real symmetric matrices. The copositive cone $\mc{C}_n$ is defined as:
\begin{align}\label{eq: def_cop}
\mc{C}_n=\left\{X\in\mc{S}_n\;|\; \mathbf{y}^\top X \mathbf{y} \geq0\textrm{ for all } \mathbf{y}\in\mathbb{R}^n_+\right\}.
\end{align}
The dual cone of $\mc{C}_n$ is the completely positive cone $\mc{C}_n^*$:
\begin{align}\label{eq: def_cpp}
\mc{C}_n^*=\left\{XX^\top \;|\; X\in\mb{R}^{n\times r}_+ \right\}.
\end{align}

In COP (CPP), we optimize a linear function of the matrix $X$ subject to linear constraints and $X\in\mc{C}_n$ ($X\in\mc{C}_n^*$).

Because COP and CPP are convex, strong duality holds if a regularity condition is satisfied, e.g., Slater's condition, which requires the feasible region to have an interior point.
\subsection{CPP Reformulation of MBQP} \label{ch: new_results}
In this section, we state the two CPP reformulations of MBQP given in \cite{burer2009copositive}. We also derive their COP duals.

Consider the MBQP:
\begin{subequations}\label{eq: bqp}
\begin{alignat}{4}
\mc{P}^\bqp: \min ~~&\x^\top Q \x + 2\c \x~~ \label{eq: bqp_obj}\\
\st ~~&\ajtop\x = b_j &&\forall j = 1, ..., m\label{eq: bqp_1}\\
&x_k\in \{0,1\} &&\forall k\in\mc{B}\label{eq: bqp_2}\\
&{\bf x} \in \mb{R}^n_+ \label{eq: bqp_3}
\end{alignat}
\end{subequations}
where $\mc{B}\subseteq\{1,...,n\}$ is the set of indices of the binary elements of $\x$. Without loss of generality we assume the matrix $Q$ is symmetric.

\cite{burer2009copositive} gives two CPP reformulations of $\mc{P}^\bqp$, (C) and (C'), which here we refer to as $\mc{P}_o^\cpp$ and  $\mc{P}_r^\cpp$. $\mc{P}_o^\cpp$ is obtained by squaring the linear constraints and substituting lifted variables for the bilinear terms.
\begin{subequations}\label{eq: cpp_origin}
\begin{alignat}{4}
\mc{P}_o^\cpp : \min ~~&\tr(Q X) + 2\c {\bf x}\label{eq: cpp_origin_obj}\\
\st ~~& \ajtop{\bf x} = b_j &&\forall j = 1, ..., m\label{eq: cpp_origin_1}\\
& \ajtop X \aj = b_j^{2} &&\forall j = 1, ..., m.\label{eq: cpp_origin_2}\\
& x_k = X_{kk} &&\forall k\in\mc{B}\label{eq: cpp_origin_3}\\
&\left[\begin{array}{cc}
1 & {\bf x}^\top \\
{\bf x} & X
\end{array}\right]\in\mc{C}_{n+1}^*.\label{eq: cpp_origin_5}
\end{alignat}
\end{subequations}

We now derive the dual of \eqref{eq: cpp_origin}. For convenience, define:
\begin{align}\label{eq: define_vars}
Y = \left[\begin{array}{cc}
1 & {\bf x}^\top \\
{\bf x} & X
\end{array}\right]; \
\tilde{Q}=\left[
\begin{array}{cc}
0&  \bm{0}^{1\times n} \\
\bm{0}^{n\times 1}& Q
\end{array}
\right]; \
C=\left[
\begin{array}{cc}
0& \c \\
\mathbf{c} & \bm{0}^{n\times n}
\end{array}
\right]; \
A_j = \left[
\begin{array}{cc}
0& 1/2 \ajtop \\
1/2 \aj & \bm{0}^{n\times n}
\end{array}
\right],\ \forall j=1,...,m.
\end{align}

We also define the matrices $\tilde{A}_j = [0, \ajtop]^\top[0, \ajtop],\;\forall j = 1,...,m$, and, for each $k\in\mc{B}$, $B_k$ such that
\begin{align*}
(B_k)_{l_1,l_2}=\left\{\begin{array}{ll}
1/2 & \textrm{if } l_1=k+1,\; l_2=1~\text{or}~l_1=1,\; l_2=k+1\\
-1 & \textrm{if } l_1=k+1,\; l_2=k+1\\
0 & \textrm{otherwise}
\end{array}\right..
\end{align*}
Then, $\mc{P}_o^\cpp$ can be written as
\begin{subequations}\label{eq: cpp_origin_lifted}
\begin{alignat}{4}
\min \quad& \textrm{Tr}(\tilde{Q} Y) + \textrm{Tr}(C Y)\label{eq: cpp_origin_lifted_obj}\\
\st \quad& \textrm{Tr}(A_j Y)=b_j\quad \forall j=1,...,m\label{eq: cpp_origin_lifted_1}\\
&\textrm{Tr}(\tilde{A}_j Y)=b_j^2 \quad \forall j=1,...,m\label{eq: cpp_origin_lifted_2}\\
&\textrm{Tr}(B_k Y)=0\quad \forall k\in\mc{B}\label{eq: cpp_origin_lifted_3}\\
&Y\in\mc{C}_{n+1}^*.\label{eq: cpp_origin_lifted_4}
\end{alignat} 
\end{subequations}

Let $\pmb{\gamma}^o$, $\pmb{\beta}^o$, $\pmb{\delta}^o$,and $\Omega^o$ be the respective dual variables of constraints \eqref{eq: cpp_origin_lifted_1} - \eqref{eq: cpp_origin_lifted_4}. The dual of $\mc{P}_o^\cpp$ is the following COP:
\begin{subequations}\label{eq: cop_origin}
\begin{alignat}{4}
\mc{P}_o^\cop : \max_{\pmb\gamma^o,\pmb\beta^o,\pmb\delta^o,\Omega^o} \quad&\sum_{j=1}^m\left(\gamma^o_j b_j+\beta^o_j b_j^2\right) \label{eq: cop_origin_obj}\\
\st \quad & \tilde{Q} + C - \sum_{j = 1}^m \gamma^o_j A_j - \sum_{j = 1}^m \beta^o_j \tilde{A}_j - \sum_{k\in\mc{B}}\delta^o_k B_k - \Omega^o = 0\label{eq: cop_origin_2}\\
& \Omega^o \in\mc{C}_{n+1}.\label{eq: cop_origin_3}
\end{alignat}
\end{subequations}

\cite{burer2009copositive} shows that $\mc{P}^\bqp$ and $\mc{P}_o^\cpp$ are equivalent, in the sense that (i) $\opt(\mc{P}^\bqp) = \opt(\mc{P}_o^\cpp)$, where $\opt(\cdot)$ is the optimal objective value; and (ii) if $(\xstar, X^*)$ is optimal for $\mc{P}_o^\cpp$, then $\mathbf{x}^*$ is in the convex hull of optimal solutions of $\mc{P}^\bqp$. The second point indicates that $\mathbf{x}^*$ is not necessarily feasible for $\mc{P}^\bqp$, i.e., it is possible for some $x_k^*$ with $k\in\mc{B}$ to be fractional.

\begin{remark}\label{th: qip_sol}
If $\xstar$ is an optimal solution of $\mc{P}^\bqp$, then $(\xstar, \xstar \mathbf{x}^{*\top})$ is optimal for $\mc{P}_o^\cpp$.
\end{remark}

\begin{remark}\label{th: feas_bqp_x}
Let $(\xstar, X^*)$ be an optimal solution for $\mc{P}_o^\cpp$. If $Q\succeq 0$ and $\xstar$ is feasible for $\mc{P}^\bqp$, then $\xstar$ is an optimal solution of $\mc{P}^\bqp$. Note that the condition $Q\succeq 0$ ensures that the objective value of $\mc{P}^\bqp$ is optimal at $\xstar$.
\end{remark}

In Section \ref{ch: app_game} we also make use of the second CPP reformulation given in \cite{burer2009copositive}, $\mc{P}_r^\cpp$. The advantages of $\mc{P}_r^\cpp$ over $\mc{P}_o^\cpp$ are that (i) it has a lower dimensional cone constraint, $X\in\mc{C}_n^*$, which may improve computational efficiency, and (ii) it may have an interior, whereas $\mc{P}_o^\cpp$ never does. $\mc{P}_r^\cpp$ is a valid reformulation of $\mc{P}^\bqp$ if the following Key Property is satisfied:
\begin{align}\label{a: alpha}
\tag{KP} \exists\; \mathbf{y}\in\mb{R}^{m}~~\st \sum^m_{j=1}y_j \aj \geq \mathbf{0},\; \sum^{m}_{j=1} y_j b_j = 1.
\end{align}
Setting ${\pmb\alpha} := \sum^m_{j=1}y_j\aj $, the CPP is given by:
\begin{subequations}\label{eq: cpp_reform}
\begin{alignat}{4}
\mc{P}_r^\cpp : \min ~~& \tr(Q X) + 2\c \x \label{eq: cpp_obj}\\
\st ~~& \ajtop{\bf x} = b_j &&\forall j = 1, ..., m\label{eq: cpp_1}\\
& \ajtop X \aj = b_j^{2} &&\forall j = 1, ..., m\label{eq: cpp_2}\\
& x_k = X_{kk} &&\forall k\in\mc{B}\label{eq: cpp_3}\\
&x_l = \tr\left(\frac{\pmb\alpha \eltop + \el \pmb\alpha^\top}{2} X\right)&~~& \forall l = 1,...,n\label{eq: cpp_4}\\
& X\in\mc{C}_n^*.\label{eq: cpp_5}
\end{alignat}
\end{subequations}

Note that $\mc{P}_r^\cpp$, as written above, looks slightly different from (C') in \cite{burer2009copositive}. This is because we have included the variable ${\mathbf x}$ via the constraint ${\mathbf x}=X\alpha$, and made several corresponding substitutions. This makes the proofs in Section \ref{ch: app_game} more straightforward. We have written the constraint ${\mathbf x}=X\alpha$ in symmetric form in \eqref{eq: cpp_4} because taking the dual of a CPP with non-symmetric linear constraints will bring further complications; we discuss the details of this in Section \ref{ec: sym_cpp} of the e-companion. Also note that Remarks \ref{th: qip_sol} and \ref{th: feas_bqp_x}  straightforwardly extend to $\mc{P}_r^\cpp$.

Condition \eqref{a: alpha} is usually not restrictive. For example, if $\mc{P}^\bqp$ contains or implies a constraint of the form $\sum_{i\in\mc{I}}a_i x_i = b$, where $a_i\geq 0$, $i\in\mc{I}$, and $b>0$, and this is the $\ell^{\text{th}}$ constraint of $\mc{P}^\bqp$, then we can construct a valid $\y$ with $1/b$ in the $\ell^{\text{th}}$ entry and 0 elsewhere.

Despite its lack of an interior, $\mc{P}_o^\cpp$ may be useful because it does not require \eqref{a: alpha} to hold. Also, strong duality may still hold for $\mc{P}_o^\cpp$, which we often observe in our numerical examples.

Finally, we write the dual of $\mc{P}_r^\cpp$. Let $\pmb{\gamma}$, $\pmb{\beta}$, $\pmb{\delta}$, $\pmb{\xi}$, and $\Omega$ be the respective dual variables of constraints \eqref{eq: cpp_1} - \eqref{eq: cpp_5}. The dual of $\mc{P}_r^\cpp$ is the COP: 
\begin{subequations}\label{eq: cop}
\begin{alignat}{4}
\mc{P}_r^\cop : \max_{\pmb\gamma,\pmb\beta,\pmb\delta,\pmb\xi, \Omega} \quad&\sum_{j=1}^m\left(\gamma_j b_j+\beta_j b_j^2\right) \label{eq: cop_obj}\\
\st \quad&2\mathbf{c}-\sum_{j=1}^m \gamma_j \aj -\sum_{k\in\mc{B}}\delta_k \ek - \sum_{l=1}^n \xi_l \el = {\mathbf 0}\label{eq: cop_1} \\
& Q - \sum_{j = 1}^m \beta_j \aj \ajtop + \sum_{k\in\mc{B}}\delta_k \ek \ektop + \sum^n_{l=1} \xi_l \frac{\pmb\alpha \eltop + \el \pmb\alpha^\top}{2} - \Omega = 0\label{eq: cop_2}\\
& \Omega\in\mc{C}_n. \label{eq: cop_3}
\end{alignat}
\end{subequations}
\section{Cutting plane algorithm}\label{ch: cut_plane}
To make use of copositive duality, we must be able to solve COPs. At present, no industrial solver can handle COPs. In the literature they are often approximately solved via SDPs (see Section \ref{ec: sdp_approx} of the e-companion). To solve COPs exactly, we design a novel cutting plane algorithm, which returns the optimal solution when it terminates. 

The cutting plane algorithm is applicable for the following general type of COPs with mixed-integer linear constraints over a copositive cone:
\begin{subequations}\label{eq: general_cop}
\begin{alignat}{4}
\max~~ &\mathbf{q}^{\top} \pmb\lambda + \tr(H^{\top} \Omega)\\
\st ~~&\mathbf{d}^\top \pmb\lambda + \tr(D_i^{\top} \Omega) = g_i &~~~&\forall i = 1,...,m\\
&\pmb\lambda\geq \boldsymbol{0}\\
&\Omega\in \mc{C}_{n_c}\label{eq: general_cop_4}\\
& \lambda_k \in \mb{Z} &~~~& \forall k\in\mc{L}\label{eq: general_cop_5}
\end{alignat}
\end{subequations}
where $\pmb\lambda$ is an $n_l$-dimensional vector, $\mc{L}$ is the index set for integer variables in $\pmb\lambda$, $\Omega\in\mb{R}^{n_c \times n_c}$. $\mc{C}_{n_c}$ is an $n_c$-dimensional copositive cone. Note that the COP problems $\mc{P}_o^\cop$ \eqref{eq: cop_origin} and $\mc{P}_r^\cop$ \eqref{eq: cop} are special cases of the COP \eqref{eq: general_cop}. We make use of the integer constraint, \eqref{eq: general_cop_5}, in solving for PNE of IPGs in Section \ref{ch: app_game}.

The algorithm starts by removing the conic constraint \eqref{eq: general_cop_4} to obtain the initial master problem, which is iteratively refined by the addition of cuts. At each iteration, we solve the master problem to obtain an optimal solution, denoted by $(\bar{\pmb\lambda}, \bar{\Omega})$. To determine whether $\bar{\Omega}$ is copositive, we employ the MIP problem proposed by \cite{anstreicher2020testing}, which checks copositivity using a recent characterization of copositive matrices given in \cite{dickinson2019new}. This MIP problem, which serves as the separation problem in our cutting plane algorithm, is given by:
\begin{subequations}\label{eq: separation}
\begin{alignat}{4}
\SP(\Omega): \ \max~~&w\\
\st ~~&\Omega \mathbf{z} \leq -w\mathbf{1} + M^\top (1- \mathbf{u})\label{eq: separation_b}\\
&\mathbf{1}^\top \mathbf{u}\geq q\label{eq: separation_c}\\
& w\geq 0\\
& \mathbf{0}\leq \mathbf{z}\leq \mathbf{u} \label{eq: linkingZandU}\\
& \mathbf{u}\in\{0,1\}^{n_c}
\end{alignat}
\end{subequations}
where $q = 1$ (or a larger integer, depending on problem structure), $\mathbf{1}$ is a vector of all ones, and $M \in \mb{R}^{n_c\times n_c}_{++}$ is a matrix of large numbers. By Theorem 2 of  \cite{anstreicher2020testing}, $\bar{\Omega}$ is copositive if and only if the optimal objective of $\SP(\bar{\Omega})$ is zero.

At any iteration, if the optimal value of the subproblem is zero, then the master problem solution is feasible and optimal for the COP \eqref{eq: general_cop}. Otherwise, as $\bar{\Omega}$ is not copositive, we add the following cut to the master problem:
\begin{align}\label{eq: cop_cut}
\mathbf{\bar{z}}^\top \Omega \mathbf{\bar{z}}\geq 0
\end{align}
where $\mathbf{\bar{z}}$ is an optimal solution of $\SP(\bar{\Omega}$).
\begin{proposition}\label{th: cop_cut}
If the optimal value of $\SP(\bar{\Omega})$ is nonzero, then \eqref{eq: cop_cut} cuts off $\bar{\Omega}$
\end{proposition}
Note that the cut \eqref{eq: cop_cut} does not eliminate any feasible solutions from \eqref{eq: general_cop}. This is because for any $\mathbf{\bar{z}}\in\mb{R}^{n_c}_+$, any copositive matrix $\Omega$ satisfies $\mathbf{\bar{z}}^\top \Omega \mathbf{\bar{z}}\geq 0$.

Since $\SP(\bar{\Omega})$ is an MIP, we can strengthen its LP relaxation to improve its computational efficiency. \cite{anstreicher2020testing} suggests doing so by setting $M_{ik} = 1 + \sum_{j = 1, ~j\neq i}^{n_c}\{\bar{\Omega}_{ij}: \bar{\Omega}_{ij} > 0\},\;\forall k = 1,...,n_c$. Another way to strengthen the separation problem is to let $q = 2$, which is valid if $\diag(\hat{\Omega}) \geq 0$ \citep{anstreicher2020testing}. Fortunately, we can always include the constraint $\diag(\hat{\Omega}) \geq 0$ and set $q=2$ because $\Omega\in\mc{C}_{n_c}$ in the initial master problem.

In some cases the master problem is unbounded at initialization. There are several ways to deal with this. One is to impose a large upper bound on the elements of $\Omega$. This bound can be gradually relaxed and then removed eventually.

If the original MBQP can be solved with reasonable efficiency and \eqref{eq: general_cop} is the dual of its completely positive reformulation, we can use strong duality in several ways. Suppose $x^*$ is the optimal solution of the MBQP.
\begin{itemize}
    \item If the optimal objective of the MBQP and the cutting plane algorithm are equal, then the algorithm has terminated at the optimal solution.
    \item Suppose $x^*$ is optimal for the MBQP. Then we can tighten the master problem by adding the complementary slackness constraint $\tr(x^*x^{*\top} \Omega) = 0$.
\end{itemize}
We can also tighten the master problem by incorporating the semidefinite relaxations of the copositive cone given in \citep{parrilo2000structured}.

We were unfortunately unable to prove that the cutting plane algorithm terminates in finite steps.  However, note that the other exact algorithm for COPs, the simplicial partition method \citep{bundfuss2009adaptive}, is also shown to be exact only in the limit. We find that in our numerical experiments, the cutting plane algorithm usually converges in a reasonable computational time, and when it does not the approximate solution obtained from the last iteration is often still useful.
\section{Pricing unit commitment}\label{ch: app_pricing}
In this section we use copositive duality to design a pricing mechanism for UC. UC optimally schedules the startups and shutdowns of the generators in a power system, typically over a 48 hour horizon. The problem is usually formulated as an MIP, in which the startup and shutdown decisions are binary variables. We reformulate the MIP as a CPP in Section \ref{subch: uc_reform}, and use the dual to define prices in Section \ref{subch: pricing_mechanism}, which we modify in Section \ref{subch: mod_cop_pricing} to guarantee revenue adequacy for individual generators.

Let $\mc{G}$ be the set of generators and $\mc{T}$ the set of time periods. $c^{p}_g$, $c^{u}_g$, and $d_t$ are respectively the production cost, startup cost, and the load. For generator $g\in\mc{G}$ at time $t\in\mc{T}$, $p_{gt}$ is the production level, $u_{gt}$ is the binary decision to startup, $z_{gt}$ is equal to one if online and zero if offline, and $\psi_{jgt}$ is the slack variable for the constraint that is not power balances.

UC can be written \citep{carrion2006computationally, taylor2015convex}:
\begin{subequations}\label{eq: uc1_reform}
\begin{alignat}{4}
\uc : \quad\min \quad& \sum_{g\in\mc{G}}\sum_{t\in\mc{T}} \left(c^{p}_g p_{gt} + c^{u}_g u_{gt}\right)\label{eq: uc1_reform_1}\\
\st \quad& \sum_{g\in\mc{G}}p_{gt} = d_t && \forall t\in\mc{T}\label{eq: uc1_reform_2} \\
&\aphitop_{jgt} \x = b_{jgt} && \forall j= 1,...,m, g\in\mc{G}, t\in\mc{T}\label{eq: uc1_reform_3}\\
& z_{gt} \in\{0,1\}&&\forall g\in\mc{G}, t\in\mc{T},
\end{alignat}
\end{subequations}
where ${\bf x}^\top = ({\mathbf u^\top}, {\mathbf z}^\top, {\mathbf p}^\top, \pmb \psi^\top)$, with bold letters of variables denoting vectors. For example, ${\mathbf u}$ denotes the vector of variables $u_{gt}$ for all $g\in\mc{G}, t\in\mc{T}$. The objective \eqref{eq: uc1_reform_1} is the total production and startup cost. Constraints \eqref{eq: uc1_reform_2} ensure that the total production satisfies the load at each hour. Constraints \eqref{eq: uc1_reform_3} are individual generator's operational constraints, which can include production level constraints, minimum up/down time, ramping constraints, and energy storage. The binary constraint on $u_{gt}$ is implied by constraints that link $u_{gt}$ and $z_{gt}$ in \eqref{eq: uc1_reform_3}.

\subsection{CPP Reformulation of UC}\label{subch: uc_reform}
We reformulate \eqref{eq: uc1_reform} in the form of $\mc{P}^{\text{CPP}}_o$. Let $X$ be the matrix of lifted variables for $\x$. $Y$ is as defined in \eqref{eq: define_vars}. {To make the correspondence between elements of $X$ and variables in vector $\x$ more explicit, we denote by $X^{vw}_{k,q}$ the element of $X$ corresponding to the row of the $v_k$ variable and the column of the $w_q$ variable. That is, $X^{vw}$ represents the block of $X$ with rows corresponding to the variables $v$ and columns to the variables $w$ by $X^{vw}$. 
For example, for the UC above, we have 
$$X=\ \begin{bmatrix}
X^{uu}_{11,11} & X^{uu}_{11,12} & ... & X^{u\psi}_{11,rl}\\
X^{uu}_{12,11}  & X^{uu}_{12,12} &  & \vdots \\
\vdots  &  & \ddots  & \vdots \\
X^{\psi u}_{rl,11} & .. & ... & X^{\psi\psi}_{rl,rl}
\end{bmatrix}.$$}
Let $\alambda_t$ be the coefficient vector for the left-hand side of constraints \eqref{eq: uc1_reform_2}. The CPP reformulation is as follows: 
\begin{subequations}\label{eq: uc1_cpp1}
\begin{alignat}{4}
\uc^\cpp : \min \quad& \sum_{g\in\mc{G}}\sum_{t\in\mc{T}} \left(c^{p}_g p_{gt} + c^{u}_g u_{gt}\right)\hspace{1cm}\label{eq: uc1_cpp1_0}\\
\st \quad& \sum_{g\in\mc{G}}p_{gt} = d_t && \forall t\in\mc{T} &&&(\lambda_t)\label{eq: uc1_cpp1_1}\\
&\aphitop_{jgt} \x = b_{jgt} && \forall j= 1,...,m, g\in\mc{G}, t\in\mc{T}\hspace{1cm}&&&(\phi_{jgt})\label{eq: uc1_cpp1_2}\\
&\tr(\alambda_t \alambdatop_t X) = d_t^2&&\forall t\in\mc{T}&&&(\Lambda_t)\label{eq: uc1_cpp1_7}\\
&\tr(\aphi_{jgt} \aphitop_{jgt} X) = b_{jgt}^2&&\forall j= 1,...,m, g\in\mc{G}, t\in\mc{T}&&&(\Phi_{jgt})\label{eq: uc1_cpp1_8}\\
&z_{gt} = Z_{gt} &&\forall g\in\mc{G},t\in\mc{T}&&&(\delta_{gt})\label{eq: uc1_cpp1_13}\\
&Y\in\mc{C}_{n+1}^*&& &&&(\Omega).\label{eq: uc1_cpp1_14}
\end{alignat}
\end{subequations}
Dual variables are shown to the right of the constraints. 

The dual of \eqref{eq: uc1_cpp1} is:
\begin{subequations}\label{eq: uc1_cop}
\begin{alignat}{4}
\uc^\cop: \max~~&\sum_{t\in\mc{T}}\left(d_t \lambda_t  + d_t^2 \Lambda_t  + \sum_{j = 1}^m \sum_{g\in\mc{G}}\left(b_{jgt}\phi_{jgt} + b_{jgt}^2 \Phi_{jgt} \right)\right)\label{eq: uc1_cop_obj}\\ %
\st ~~&(\pmb \lambda, \pmb \phi, \pmb \Lambda, \pmb \Phi, \pmb \delta, \Omega)\in\mc{F}^{\text{COP}},\label{eq: uc1_cop_constr}
\end{alignat}
\end{subequations}
where $\mc{F}^{\text{COP}}$ denotes the feasible region of the dual problem, which can be written in the form of constraints \eqref{eq: cop_origin_2} - \eqref{eq: cop_origin_3}.
\subsection{Copositive Duality Pricing}\label{subch: pricing_mechanism}
We now describe CDP, a pricing mechanism for UC. CDP is budget balanced and, under certain conditions, individually rational. Let $\x^*$ be the optimal solution of $\uc$, and set $X^* = \x^* \x^{*\top}$. According to Remark \ref{th: qip_sol}, $(\x^*, X^*)$ is an optimal solution to $\uc^\cpp$ \eqref{eq: uc1_cpp1}. The CDP mechanism is defined as follows:
\begin{definition}[CDP]
\label{def:CDP}
Let $(\pmb \lambda^*, \pmb \phi^*, \pmb \Lambda^*, \pmb \Phi^*)$ be an optimal solution for $\uc^\cop$. Under the CDP mechanism, at hour $t$ the system operator (SO):
\begin{itemize}
\item collects $\pi^{\text L}_t =  \lambda_t^*d_t  +  \Lambda_t^*d_t^2  +  \sum_{j=1}^m \sum_{g\in\mc{G}}\left(b_{jgt}\phi_{jgt}^* + b_{jgt}^2 \Phi_{jgt}^* \right)$ from the load, and 

\item pays $\pi^{\text G}_{gt} = \lambda_t^* p^*_{gt} + \Lambda^*_t \Xpstar + \sum_{j=1}^m\left(\phi^{*}_{jgt}\aphi_{jgt} \x^* + \Phi^{*}_{jgt}\tr(\aphi_{jgt} \aphitop_{jgt} X^*)\right) + \sum_{g\in -g}f_{g g' t}$ to generator $g$ at time $t$, where $-g = \mc{G}\sm \{g\}$. $f_{g g' t}$ is the share of $g$'s revenue from the cross-term payment $2\Lambda_t^* \Xpp$. It must satisfy $f_{g g' t} + f_{g' g t} = 2\Lambda_t^* \Xpp$, and if $\Xpp = 0$, then $f_{g g' t} = 0$. 
\end{itemize}
\end{definition}

$\pi^{\text L}_t$ consists of volumetric price payments, $d_t \lambda_t^*$  and $d_t^2 \Lambda_t^*$, and payments that depend on the shadow prices of operational constraints with non-zero right-hand sides. Note that the quadratic term $d_t^2 \Lambda_t^*$ corresponds to the lifted power balance \eqref{eq: uc1_cpp1_7}. The shadow price payments for operational constraints are roughly comparable to transmission congestion rents. They could represent, for example, a payment corresponding to a ramping constraint, which, if loosened, would improve the objective. 


$\pi^{\text G}_{gt}$ depends on generators' optimal production levels and on/off statuses. It is obtained by summing the products of the left-hand sides of constraints \eqref{eq: uc1_cpp1_1} - \eqref{eq: uc1_cpp1_8} with their corresponding dual prices. As in the load payment, $\lambda_t^* p^*_{gt}$ and $\Lambda^*_t \Xpstar$ are volumetric payments for which the prices are uniform for all generators. $\Lambda^*_t \Xpstar$ corresponds to the lifted power balance \eqref{eq: uc1_cpp1_7}. The terms in the first sum similarly correspond to shadow prices of various constraints. The shadow prices are not uniform as they depend on the generator index, $g$. The second sum corresponds to the off-diagonal entries in the lifted matrix $X^*$. Since this payment involves two generators, it is an open question as to how this payment should be divided between those generators. For example, we can assign it to the generator that loses money, or divide it evenly as in Section \ref{subch: mod_cop_pricing}.

The following example examines non-uniform payments resulting from shadow prices.  We use this example in our numerical experiments in Section \ref{subsec: electricity_markets_results}.

\begin{example}\label{case: simple_uc}
Suppose constraint \eqref{eq: uc1_reform_3} is given by:
\begin{subequations}\label{eq: specific_uc}
\begin{alignat}{4}
& u_{gt} -\theta_{gt} = z_{gt} - z_{g, t-1} &~~~& \forall g\in\mc{G}, t\in\mc{T}\sm\{1\}\label{eq: specific_uc_1}\\
& p_{gt} - \mu_{gt}= z_{gt} \underline{P}_g && \forall g\in\mc{G}, t\in\mc{T}\label{eq: specific_uc_2}\\
& p_{gt} + \gamma_{gt} = z_{gt} \overline{P}_g && \forall g\in\mc{G}, t\in\mc{T}\label{eq: specific_uc_3}\\
& z_{gt} + \eta_{gt} = 1 && \forall g\in\mc{G}, t\in\mc{T}, \label{eq: specific_uc_4}
\end{alignat}
\end{subequations}
where $\theta_{gt}, \mu_{gt}, \gamma_{gt}$, and $\eta_{gt}$ are slack variables, and $\underline{P}_g$ and $\overline{P}_g$ are respectively the lower and upper production limits of generator $g\in\mc{G}$. Constraint \eqref{eq: specific_uc_4} implies upper bounds on $z_{gt}$, and is the only constraint with a non-zero constant on the right-hand-side. 

The portion of the payment resulting from constraint \eqref{eq: specific_uc_4} and the corresponding lifted constraint are included in $\pi^{\text G}_t$. We refer to the shadow prices of these constraints the {\it availability prices}, as they signal the availability of the generators at each hour. These are the only non-uniform prices in this example. 
\end{example}

The following theorem shows that CDP is revenue neutral for both the SO and the generators. Note that the individual generators are in general not revenue neutral under CDP.
\begin{theorem}\label{th: rev_adequacy}
If strong duality holds for $\uc^\cpp$, then CDP balances the revenue and the aggregate cost of the generators. It also balances the revenue and payment of the SO.
\end{theorem}
Next, we provide a sufficient condition guaranteeing the individual rationality of CDP. A pricing mechanism is individually rational if each market participant has no incentive to deviate from the centrally optimal solution.

The profit-maximization problem of generator $g$ is given by:
\begin{subequations}\label{eq: uc1_individual}
\begin{alignat}{4}
\pi_g^\gen (\p_{-g}) := \quad\max_{\xg} \quad& \sum_{t\in\mc{T}} \Big(\lambda^*_t p_{gt} + \Lambda^*_t p_{gt}^{2} + \sum_{j=1}^m\left(\phi^*_{jgt}b_{jgt} + \Phi^*_{jgt}b_{jgt}^2 \right) +\sum_{g'\in -g}f_{g g' t}\nonumber\hspace{-3cm}\\&- c^p_g p_{gt} - c^u_g u_{gt}\Big)\label{eq: uc1_individual_1}\\
\st \quad&\aphitop_{jgt} \x = b_{jgt}  && \forall j = 1,...m, t\in\mc{T}&&~~&\label{eq: uc1_individual_2}
\end{alignat}
\end{subequations}
where $\p_{-g}$ denotes the production decision of all generators other than $g$. $\xg\in\mb{R}^{n_g}$ denotes the portions of $\x$ that corresponds to generator $g\in\mc{G}$. The first five terms in the objective represent the total revenue of generator $g$. $\pi_g^\gen (\p_{-g})$ should not contain entries of $X$, which are lifted variables that do not have direct physical meanings. As discussed in Remark \ref{th: qip_sol}, at optimality we can equivalently use $(p^{*}_{gt})^2$ in place of $\Xpstar$ and $p^*_{gt} p^*_{g' t}$ in place of $\Xppstar$.

If the optimal solution of the profit-maximization problem \eqref{eq: uc1_individual} matches the corresponding portion of the optimal solution of $~\uc$ for all generators $g\in\mc{G}$, then the market mechanism is individually rational.

Theorem \ref{th: ind_rational} provides a sufficient condition for when CDP is individually rational.
\begin{theorem}\label{th: ind_rational}
Assume strong duality holds for $\uc^\cpp$, and let $(\x^*, X^*)$ be an optimal solution. If all items in $X^*$ in the forms of $\Xvwstar, \forall { v,w}, g_1 \neq g_2, g_1,g_2\in\mc{G}, t\in\mc{T}$ are equal to zero, then the market mechanism with CDP is individually rational.
\end{theorem}

The condition $\Xvwstar = 0, \forall {v,w}, g_1 \neq g_2, g_1,g_2\in\mc{G}, t\in\mc{T}$ ensures that the conic constraint in $\uc^\cpp$ is decomposable by $g$, which is a key property needed in the proof.

More generally, there may be groups of generators that are coupled to each other, but not to those of other groups. We refer to this as subset rationality. There could also, for example, be groups of time periods that are similarly decoupled. It may be possible to identify more complicated decomposable structures using matrix completion \citep{drew1998completely}. This is a topic of future work.

\subsection{Ensuring Individual Revenue Adequacy}\label{subch: mod_cop_pricing}
CDP does not guarantee each individual generator's revenue adequacy, i.e., nonnegative profit. This is also the case with some other schemes, including RP and CHP. In CDP, we can enforce revenue adequacy by adding constraints directly to the dual of $\uc^\cpp$.

Revenue adequacy can be enforced through the non-uniform prices $\phi_{jgt}^*$ and $\Phi_{jgt}^*$, which can be different for each generator, and/or through the uniform prices $\lambda_t$ and $\Lambda_t$. 
We use uniform pricing because it is easier to implement in practice. In Section \ref{ec: ind_rev_avail_price} of the e-companion, we present a different version that uses availability prices as well.

We enforce revenue adequacy by adding a new constraint to the dual problem, which ensures that the payment to each generator is no less than its costs. The resulting augmented dual problem is given by:
\begin{subequations}\label{eq: uc1_cop_profitable}
\begin{alignat}{4}
\max~~&\sum_{t\in\mc{T}}\left(d_t \lambda_t + d_t^2 \Lambda_t + \sum_{j=1}^m \sum_{g\in\mc{G}}( b_{jgt}\phi_{jgt} +  b_{jgt}^2 \Phi_{jgt})\right)\label{eq: uc1_cop_profitable_obj}\\ 
\st ~~&\sum_{t\in\mc{T}} \left(p^*_{gt}\lambda_t + p^{*2}_{gt} \Lambda_t + \sum_{g'\in -g} p_{gt}^*p_{g't}^*\Lambda_t \right) \geq \sum_{t\in\mc{T}}\left(c^p_g p^*_{gt} + c^u_g u^*_{gt}\right) &~~~&\forall g\in\mc{G} \label{eq: uc1_cop_profitable_1}\\
&(\pmb \lambda, \pmb \phi, \pmb \Lambda, \pmb \Phi, \pmb \delta, \Omega)\in\mc{F}^{\text{COP}}.\hspace{-1cm}&&\label{eq: uc1_cop_profitable_constr}
\end{alignat}
\end{subequations}
The objective \eqref{eq: uc1_cop_profitable_obj} and the constraint \eqref{eq: uc1_cop_profitable_constr} are the same as in the original dual problem, $\uc^\cop$. The left-hand side of constraint \eqref{eq: uc1_cop_profitable_1} is the total revenue of generator $g$, assuming that the cross-term payments are divided evenly between generators. The right-hand side of \eqref{eq: uc1_cop_profitable_1} is the total cost of generator $g$.

The prices are computed by solving \eqref{eq: uc1_cop_profitable}. The optimal values of primal variables, $p_{gt}^*$ and $u_{gt}^*$, could be obtained by solving $\uc$. We refer to this pricing mechanism as Revenue-adequate CDP (RCDP). More formally, the RCDP mechanism is defined as follows.
\begin{definition}[RCDP Mechanism]
Let $(\pmb\lambda^*, \pmb\Lambda^*)$ be an optimal solution for \eqref{eq: uc1_cop_profitable}. Under the RCDP mechanism, at hour $t$ the SO:
\begin{itemize}
\item collects $\pi^{{\text L}'}_t = \lambda_t^*d_t   +  \Lambda_t^*d_t^2$ from the load, and 

\item pays $\pi^{{\text G}'}_{gt} = \left(\lambda_t^* p^*_{gt} + \Lambda^*_t \Xpstar \right) + \sum_{g\in-g} p_{gt}^*p_{g't}^*\Lambda_t$ to the generator $g$. 
\end{itemize}
\end{definition}

A benefit of RCDP is that all generators are paid the same price in each hour, so that more efficient generators with lower costs earn higher profit. The load also pays on a purely volumetric basis. The SO remains revenue neutral, in that $\pi^{{\text L}'}_t=\pi^{{\text G}'}_t$. This can be proved by adapting the proof of Theorem \ref{th: rev_adequacy}.

Because the addition of constraints \eqref{eq: uc1_cop_profitable_1} restrict the original dual problem \eqref{eq: uc1_cop}, it is worth checking if the problem is still feasible. We known that if either the primal or dual is feasible, bounded, and has an interior point, then the other is also feasible \citep{luenberger2015linear}. We therefore focus on the dual of \eqref{eq: uc1_cop_profitable}. Let $q_g \geq 0$ be the dual multiplier of \eqref{eq: uc1_cop_profitable_1}. The dual of \eqref{eq: uc1_cop_profitable} is the following:
\begin{subequations}\label{eq: uc1_cpp_profitable}
\begin{alignat}{4}
\min \quad& \sum_{g\in\mc{G}}\sum_{t\in\mc{T}} \left(c^p_g p_{gt} + c^u_g u_{gt}- (c^p_g p^*_{gt} + c^u_g u^*_{gt})q_g\right)\label{eq: uc1_cpp_profitable_0}\\
\st \quad& \sum_{g\in\mc{G}}\left(p_{gt} - p^*_{gt}q_g\right) = d_t && \forall t\in\mc{T} \label{eq: uc1_cpp_profitable_1}\\
&\tr(\alambda_t \alambdatop_t X) - \sum_{g\in\mc{G}} \left(p^{*2}_{gt} + \sum_{g'\in -g} p_{gt}^*p_{g't}^*\right)q_g= d_t^2&~~~&\forall t\in\mc{T}\label{eq: uc1_cpp_profitable_2}\\
& \eqref{eq: uc1_cpp1_2}, \eqref{eq: uc1_cpp1_8}, \eqref{eq: uc1_cpp1_13}\\
&q_{g} \geq 0 &&g\in\mc{G}\\
&Y\in\mc{C}_{n+1}^*. && \label{eq: uc1_cpp_profitable_4}
\end{alignat}
\end{subequations}

Observe that this is similar to $\uc^\cpp$, but with additional terms multiplying $q_g$ in the objective and constraints. Note that \eqref{eq: uc1_cpp_profitable} includes both $p_{gt}$, a variable, and $p_{gt}^*$, part of the solution to $\uc^\cpp$.

The CPP \eqref{eq: uc1_cpp_profitable} is feasible because we recover $\uc^\cpp$ by setting all $q_g$ to zero. It is bounded due to the equality \eqref{eq: uc1_cpp_profitable_1}, and the fact that the production level $p_{gt}$ is usually bounded by generator capacity. If the problem has an interior, then its dual, \eqref{eq: uc1_cop_profitable}, is feasible. 
\begin{remark}
Since \eqref{eq: uc1_cpp_profitable} is a relaxation of $\uc^\cpp$, \eqref{eq: uc1_cpp_profitable} should have an interior if $\uc^\cpp$ has an interior. Also, if $\uc^\cpp$ has an interior, then $\uc^\cop$ is feasible. Thus, the RCDP problem \eqref{eq: uc1_cop_profitable} is feasible whenever $\uc^\cop$ is feasible.
\end{remark}

\section{Mixed-Binary Quadratic Games}\label{ch: app_game}
In this section, we use copositive duality to characterize the equilibria of mixed-binary quadratic games (MBQGs). In an MBQG, each player solves an optimization of the form $\mc{P}^\bqp$, which are coupled through the objectives, but not the constraints. The binary decisions preclude the direct use of existence and uniqueness results for continuous games, and the use of the KKT conditions for obtaining equilibria. We surmount these issues by converting each player's MBQP to a CPP, as described in Section \ref{ch: new_results}. We give conditions for existence and uniqueness in Section \ref{subch: pne_exist}, and compute PNE using the KKT conditions in Section \ref{subch: compute_ne}.

We study an $n$-person MBQG, $\mc{G}^\bqp = \langle\mc{I}, (\mcXi)_{i\in\mc{I}}, (\xith)_{i\in\mc{I}}\rangle$, where $\mc{I}$ is the set of indices for the players and $|\mc{I}| = n$, $\mcXi$ is the strategy set of player $i$, and $\xith\in\mcXi$ is the $n_i$-dimensional decision vector of player $i$, which contains both continuous and binary variables. We assume that $\mcXi$ is bounded throughout this section. If it is not, we can always add constraints to make it bounded without changing the optimal solution. $\mcXi$ is also closed because there are no strict inequality constraints in $\mc{P}^\bqp_i(\xni)$ below, and therefore compact. Let $\mcXni$ and $\xni$ denote the strategy sets and strategies of all players except $i$. Player $i$ faces the following MBQP:
\begin{subequations}\label{eq: ip_game}
\begin{alignat}{4}
\mc{P}^\bqp_i(\xni) : \min ~~& f_i(\xith, \xni)\label{eq: ip_game_obj}\\
\st ~~& \aijtop\xith = \bij &~~&\forall j = 1, ..., m_i\label{eq: ip_game_1}\\
&x_{ik}\in \{0,1\} &&\forall k\in\mc{B}_i\label{eq: ip_game_3}\\
&\xith \in \mb{R}_+^{l}.\label{eq: ip_game_4}
\end{alignat}
\end{subequations}

$f_i(\xith, \xni) = \x^\top Q^{(i)}\x + 2 {\bf c}^{(i)\top} {\bf x}$ is the quadratic payoff function of player $i$, $\x$ the vector of all strategies, and $Q^{(i)} \in \mc{S}_{n_i}$ a symmetric matrix. Constraints \eqref{eq: ip_game_1} are linear and independent of the decisions of other players. $\mc{B}_i$ is the subset of indices corresponding to the binary decisions of player $i$.

A pure-strategy NE (PNE) of $\mc{G}^\bqp$ is a vector of actions ${\bf x}^*\in(\mcXi)_{i\in\mc{I}}$ such that for each player $i$: 

$$f_i(\xstari, \xstarni)\leq f_i(\xith, \xstarni),\;\forall \xith\in\mcXi,$$
which is equivalently stated
\begin{align}
f_i(\xstari, \xstarni) = \text{opt}(\mc{P}^\bqp_i(\xstarni)).
\end{align}
\subsection{Reformulation in terms of CPP}\label{ch: reduction}
We reformulate $\mc{G}^\bqp$ in terms of CPP, and show how the PNE of one can be obtained from the other, and vice versa. Throughout this section, we use the CPP reformulation $\mc{P}_r^\cpp$ because it has several advantages over $\mc{P}_o^\cpp$. For this reason we assume that \eqref{a: alpha} holds throughout. However, we could in principle derive all of the results in this section using $\mc{P}_o^\cpp$, in which case \eqref{a: alpha} need not hold.

We denote the reformulated game $\mc{G}^\cpp= \langle\mc{I}, (\mcXip)_{i\in\mc{I}}, (\xith, \Xith)_{i\in\mc{I}}\rangle$, where $\mcXip$ is the strategy set of player $i$. Due to the relationship between the feasible regions of $\mc{P}^\bqp$ and $\mc{P}^\cpp$, as established by Corollary 2.5 of \cite{burer2009copositive}, the compactness of $\mcXi$ implies the compactness of $\mcXip$. Player $i$ solves:
\begin{subequations}\label{eq: cpp_game}
\begin{alignat}{4}
\mc{P}^\cpp_i(\xni) : \min ~~&\tr(Q_{ii}^{(i)}\Xith) + 2\xnitop Q^{(i)}_{-i,i} \xith + \xnitop \Qi_{-i,-i} \xni + 2\ci {\bf x}\label{eq: cpp_game_obj}\\
\st ~~& \aijtop\xith = \bij &&\forall j = 1, ..., m_i\label{eq: cpp_game_1}\\
& \aijtop \Xith\aij = b_j^{(i)2} &&\forall j = 1, ..., m_i\label{eq: cpp_game_2}\\
& x_{ik} = X_{i, kk} &&\forall k\in\Bi\label{eq: cpp_game_3}\\
& x_{il} = \tr\left(\frac{\alphai \eltop + \el \alphaitop}{2} \Xith\right) && \forall l = 1,...,n_i\label{eq: cpp_game_4}\\
& \Xith\in\mc{C}_{n_i}^*. \label{eq: cpp_game_5}
\end{alignat}
\end{subequations}
Note that $\xith$ and $X_i$ are the only variables in $\mc{P}^\cpp_i(\xni)$, and $\xni$ is treated as given. $\alphai$ is defined for $\mc{P}^\bqp_i$ as in \eqref{a: alpha}.

Let $F_i(\xith, \Xith, \xni)$ denote the objective function of $\Pith^\cpp(\xni)$. A PNE of $\mc{G}^\cpp$ is a set of strategies, $(\mathbf{x}^*_i, X^*_i)_{i\in\mc{I}}\in(\mcXip)_{i\in\mc{I}}$, such that for each player $i$:
\begin{align}
F_i(\xstari, \Xstari, \xstarni) = \text{opt}(\mc{P}^{\text{CPP}}_i(\xstarni)).
\end{align}

Any game in the form of $\mc{G}^\bqp$ can be written in the form of $\mc{G}^\cpp$. Theorem \ref{th: games_equiv} relates the equilibria of the games.

\begin{theorem}\label{th: games_equiv}
Let $\mathbf{x}^*$ be a PNE of $\mc{G}^\bqp$. Then $(\xstari, \xstari \xstaritop)_{i\in\mc{I}}$ is a PNE of $\mc{G}^\cpp$. 

Conversely, if $(\xstari, \Xstari)_{i\in\mc{I}}$ is a PNE of $\mc{G}^\cpp$, and if $(\Qi_{ii}\succeq 0)_{i\in\mc{I}}$, $(\xstarik\in\{0,1\},\;\forall k\in\mc{B}_i)_{i\in\mc{I}}$, then $\mathbf{x}^*$ is a PNE of $\mc{G}^\bqp$. 
\end{theorem}
{\it Remark. }The requirement $(\xstarik\in\{0,1\},\;\forall k\in\mc{B}_i)_{i\in\mc{I}}$ ensures that a PNE of $\mc{G}^\cpp$ can be transformed into $\mc{G}^\bqp$. We can replace this requirement with the (possibly more restrictive) condition that the matrix $\begin{bmatrix}
1 & \x^{*\top}_i\\
\xstari & \Xstari
\end{bmatrix}$ is rank one for all $i\in\mc{I}$. This is similar to the condition given in \cite{ahmadi2020semidefinite} to ensure the validity of an SDP reformulation of a 2-person bimatrix game. 

The second part of Theorem \ref{th: games_equiv} can also be replaced with the following (more restrictive) proposition:
\begin{proposition}\label{th: games_equiv2}
If $(\xstari, \Xstari)_{i\in\mc{I}}$ is a PNE of $\mc{G}^\cpp$, and if $(\Qi_{ii}\succeq 0)_{i\in\mc{I}}$, and each one of the optimization problems $(\Pith^\bqp(\xstarni))_{i\in\mc{I}}$ has a unique optimal solution, then $\mathbf{x}^*$ is a PNE of $\mc{G}^\bqp$.
\end{proposition}

This proposition is true because if $\Pith^\bqp(\xstarni)$ has a unique optimal solution, then $\xstari$ must be feasible for $\Pith^\bqp(\xstarni)$. Thus we can follow a similar proof as in Theorem \ref{th: games_equiv} to show that $\mathbf{x}^*$ is a PNE of $\mc{G}^\bqp$. 
\subsection{Existence and Uniqueness of PNE in $\mc{G}^\bqp$}\label{subch: pne_exist}
We first prove that PNE for $\mc{G}^\bqp$ exists under some conventional assumptions. 
To prove the existence of PNE in $\mc{G}^\bqp$, first note that its corresponding CP game $\mc{G}^\cpp$ always has at least one PNE. This is because according to the classical PNE existence condition for games with convex strategy sets \citep{debreu1952social, glicksberg1952further, fan1952fixed}, $\mc{G}^\cpp$ should have a PNE as $\mc{P}^\bqp_i(\xni)$ has a convex and compact feasible region and a linear objective function.

Now we provide a sufficient condition for the existence of PNE in $\mc{G}^\bqp$.
\begin{proposition}[Existence of PNE]\label{th: exist_qip}$\mc{G}^\bqp$ has at least one PNE if
\BI
\I[(i)] Property \eqref{a: alpha} is satisfied for $\Pith^{\bqp}(\xni)$,

\I[(ii)] $\Qi_{ii}\succeq 0$, and

\I[(iii)] there exists a PNE $(\xstari, \Xstari)_{i\in\mc{I}}$ for the corresponding $\mc{G}^\cpp$ that satisfies $(\xstarik\in\{0,1\},\;\forall k\in\mc{B}_i)_{i\in\mc{I}}$.
\EI
\end{proposition}

Conditions (ii) and (iii) are necessary to ensure that PNE of $\mc{G}^\cpp$ can be transformed to the PNE of $\mc{G}^\bqp$ (see Theorem \ref{th: games_equiv}). If condition (iii) is not satisfied, then $\mc{G}^\bqp$ does not have any PNE. This is because if $\mc{G}^\bqp$ did have a PNE, then, following Theorem \ref{th: games_equiv}, we could construct a PNE for $\mc{G}^\cpp$ that does satisfy (iii).

Using Proposition \ref{th: games_equiv2}, condition (iii) in Proposition \ref{th: exist_qip} can also be replaced with the requirement that the optimization problems in $\mc{G}^\bqp$ all have unique optimal solutions, since this requirement also ensures that a PNE of $\mc{G}^\cpp$ can be converted to a PNE of $\mc{G}^\bqp$. It is worth noting that a similar requirement is proposed for the existence of PNE in two-person discrete games, which is a special type of MBQ game, by \cite{mallick2011existence}. This paper proves that for a two-person discrete game, if both players have unique best responses and a condition called Minimal Acyclicity is satisfied, then the PNE must exist.

We note that in principle the uniqueness theorem of \cite{rosen1965existence} can be applied to $\mc{G}^\cpp$. In particular, it implies that if equilibria exist, certain regularity conditions hold for the player problems to have strong duality, and the payoff functions are diagonally strictly convex, then any PNE of $\mc{G}^\bqp$ is unique. It may be possible to use this perspective to identify new classes of discrete games that have a unique equilibrium.

\subsection{Computing PNE}\label{subch: compute_ne}
A PNE of $\mc{G}^\cpp$ satisfies the KKT conditions when a certain constraint qualification, such as Slater's condition, is met. We write the KKT conditions for $\mc{G}^\cpp$ explicitly. For all $i\in\mc{I}$, we have:
\begin{subequations}\label{eq: cpgame_kkt}
\begin{alignat}{4}
&\eqref{eq: cpp_game_1} - \eqref{eq: cpp_game_5} \label{eq: cpgame_kkt_1}\\
&2Q_{-i,i}^{(i)\top}\xni + 2\mathbf{c}_i-\sum_{j=1}^{m_i} \gamma_{ij} \mathbf{a}^{(i)}_j -\sum_{k\in\mc{B}_i}\delta_{ik} \ek - \sum_{l=1}^{n_i} \xi_{il} \el = 0\label{eq: cpgame_kkt_2} \\
& Q_{ii}^{(i)} - \sum_{j = 1}^{m_i} \beta_{ij} \mathbf{a}_j^{(i)} \a^{(i)\top}_j + \sum_{k\in\mc{B}_i}\delta_{ik} \ek \ektop + \sum_{l=1}^{n_i} \xi_{il} \frac{\alphai \eltop + \el \alphaitop}{2} - \Omega_i = 0 \label{eq: cpgame_kkt_3}\\
& \Omega_i\in\mc{C}_{n_i} \label{eq: cpgame_kkt_4}\\
& \tr(\Omega_i \Xith) = 0 \label{eq: cpgame_kkt_5}
\end{alignat}
\end{subequations}
where \eqref{eq: cpgame_kkt_1} are constraints of $\mc{P}^\cpp_i(\xni)$. $\pmb{\gamma}_i$, $\pmb{\beta}_i$, $\pmb{\delta}_i$, $\pmb{\xi}_i$, and $\Omega_i$ are the respective dual variables for \eqref{eq: cpp_game_1} - \eqref{eq: cpp_game_5}. Let $\mc{L}(\xith, \Xith, \xni)$ be the Lagrangian function of $\mc{P}^\cpp_i(\xni)$, i.e.,  
\begin{align*}
\mc{L}(\xith, \Xith, \xni) =  ~&F_i(\xith, \Xith, \xni) + \sum_{i=1}^{m_i} \left( \gamma_{ij}(\bij - \aijtop\xith) + \beta_{ij} ((b_j^{(i)})^2 - \aijtop \Xith \aij)\right)\\
& + \sum_{k\in\mc{B}_i}\delta_{ik}\left(X_{i, kk} - x_{ik}\right) + \sum_{l=1}^{n_i} \xi_{il} \left(\tr\left(\frac{\alphai \eltop + \eltop \alphai}{2} X_{i}\right) - x_{il}\right) + \tr(\Omega_i \Xith).
\end{align*}
Then $\partial\mc{L}(\xith, \Xith, \xni)/\partial \xith = 0$ and $\partial\mc{L}(\xith, \Xith, \xni) / \partial \Xith= 0$ respectively correspond to constraints \eqref{eq: cpgame_kkt_2} and \eqref{eq: cpgame_kkt_3}. Constraint \eqref{eq: cpgame_kkt_4} is the dual cone constraint. Finally, constraint \eqref{eq: cpgame_kkt_5} enforces complementary slackness. 

Our goal is to obtain the PNE for $\mc{G}^\bqp$. By Proposition \ref{th: exist_qip}, if $Q^{(i)}_{ii} \succeq 0,\;\forall i\in\mc{I}$, then we can add the following extra constraints to the KKT conditions of $\mc{G}^\cpp$:
\begin{align}\label{eq: kkt_binary}
x_{ik} \in\{0,1\},\;\forall k\in\mc{B}_i, i\in\mc{I},
\end{align}
to ensure that the solution is a PNE for $\mc{G}^\bqp$. 

The KKT conditions \eqref{eq: cpgame_kkt} contain conic constraints $\Xith\in\mc{C}_{n_i}^*$ and $\Omega_i\in\mc{C}_{n_i}$, as well as a bilinear constraint $\tr(\Omega_i \Xith) = 0$. As a result, it is difficult to solve \eqref{eq: cpgame_kkt} with current technology. However, in the special case where $|B_i| = l$ in $\mc{P}^\bqp_i(\xni)$ \eqref{eq: ip_game} $\forall i \in\mc{I}$, i.e., all variables in $\mc{P}^\bqp_i(\xni)$ are binary (such as in bimatrix games), then the KKT conditions \eqref{eq: cpgame_kkt} could be replaced with a MIP problem over copositive cone constraints. We now reformulate the KKT conditions for this special case.

From the primal problem constraints \eqref{eq: cpgame_kkt_1}, we drop the completely positive conic constraint \eqref{eq: cpp_game_5}. We keep constraints \eqref{eq: cpp_game_1}, and add constraints \eqref{eq: kkt_binary}. Observe that constraints \eqref{eq: cpp_game_2} - \eqref{eq: cpp_game_5} are a relaxation of the constraint $\Xith = \xith \xitop$, or, equivalently, $X_{i,j_1, j_2} = x_{i j_1} x_{i, j_2}$. Because all variables in $\mc{P}^\bqp_i(\xni)$ are binary, we can therefore replace \eqref{eq: cpp_game_2} - \eqref{eq: cpp_game_5} with the following linear constraints $(\forall i\in\mc{I})$:
\begin{subequations}\label{eq: kkt_mip_comp1}
\begin{alignat}{4}
& X_{i, jj} = x_{ij} && \forall j = 1,...,l \\
& X_{i, j_1 j_2} \leq x_{i j_1} &&\forall j_1\neq j_2; j_1,j_2 = 1,...l\label{eq: kkt_mip_comp1_2} \\
& X_{i, j_1 j_2} \leq x_{i j_2}&&\forall j_1\neq j_2; j_1,j_2 = 1,...l \\
& X_{i, j_1 j_2} \geq x_{i j_1} + x_{i j_2} - 1&~~&\forall j_1\neq j_2; j_1,j_2 = 1,...l \\
& X_{i, j_1 j_2} \geq 0 &~~&\forall j_1\neq j_2; j_1,j_2 = 1,...l.\label{eq: kkt_mip_comp1_5}
\end{alignat}
\end{subequations}
Constraints \eqref{eq: kkt_mip_comp1_2}-\eqref{eq: kkt_mip_comp1_5} are the McCormick relaxation of the bilinear expression $X_{i,j_1, j_2} = x_{i j_1} x_{i, j_2} ~(j_1\neq j_2)$.

On the dual side, we keep all constraints \eqref{eq: cpgame_kkt_2} - \eqref{eq: cpgame_kkt_4}. We linearize the bilinear complementary slackness constraint, \eqref{eq: cpgame_kkt_5}, by defining a new matrix, $Z_i \in\mc{R}^{l\times l}$, and letting $Z_{i, j_1 j_2} = \Omega_{i, j_1 j_2} X_{i, j_1 j_2}$. Constraint \eqref{eq: cpgame_kkt_5} can then be replaced with the disjunctive:
\begin{subequations}\label{eq: kkt_mip_comp2}
\begin{alignat}{4}
&\tr(Z) = 0\\
&-m X_{i, j_1 j_2} \leq Z_{i, j_1 j_2} \leq m X_{i, j_1 j_2} &&\forall j_1,j_2 = 1,...,l\\
&\Omega_{i, j_1 j_2} - m (1 - X_{i, j_1 j_2}) \leq Z_{i, j_1 j_2} \leq \Omega_{i, j_1 j_2} + m (1 - X_{i, j_1 j_2})&~~&\forall j_1,j_2 = 1,...,l
\end{alignat}
\end{subequations}
where $m$ is a sufficiently large number. Here we use the fact that elements in $\Xith$ are binary, which is due to constraints \eqref{eq: kkt_binary} and \eqref{eq: kkt_mip_comp1}.

We can now state the reformulation of the KKT conditions \eqref{eq: cpgame_kkt} for the case where all variables are binary. We have, $\forall i\in\mc{I}$:
\begin{subequations}\label{eq: kkt_mip}
\begin{alignat}{4}
&  \aijtop \xith = b_j^{(i)} &~~&\forall j = 1, ..., m_i\\
&\eqref{eq: cpgame_kkt_2}, \eqref{eq: cpgame_kkt_3}, \eqref{eq: kkt_binary}, \eqref{eq: kkt_mip_comp1}, \eqref{eq: kkt_mip_comp2}\\
& \Omega_i\in\mc{C}_{n_i}.
\end{alignat}
\end{subequations}
This is an MIP over the copositive cone, and can be solved exactly with the cutting plane algorithm in Section \ref{ch: cut_plane} using an MIP solver. If a solution to \eqref{eq: kkt_mip} exists, then it is a PNE for $\mc{G}^\bqp$. {Note that we cannot approximate \eqref{eq: kkt_mip} with SDP solvers such as MOSEK because they cannot readily solve mixed-integer SDPs.}

\section{Numerical Results}\label{ch: results}
In this section we present our numerical experiments. Section \ref{ec: scarf} implements several pricing schemes on the Scarf's example. Section \ref{subsec: electricity_markets_results} compares different pricing schemes for a nonconvex electricity market. In Section \ref{subsec: bimatrix}, we compute the PNE of bimatrix games using their KKT conditions. We also include a comparison of our COP cutting plane algorithm with other COP algorithms in e-companion Section \ref{ec: max_clique}.

All experiments are implemented in Julia v1.0.5 using the optimization package JuMP.jl v0.20.1. The COP cutting plane algorithm was implemented using CPLEX 12.8. We use Mosek 9.1 to solve SDPs approximations in Sections \ref{ec: scarf}, \ref{subsec: electricity_markets_results}, and \ref{ec: max_clique}.

\subsection{Pricing in Scarf's Example}\label{ec: scarf}
Scarf's example is often used to compare pricing schemes for nonconvex markets. We use the modified version from \cite{hogan2003minimum} to compare CDP with RP and CHP, which are currently used by utilities in the U.S.. In the modified Scarf's example, there are three types of generators: smokestack, high technology, and medium technology. Let $\mc{G}_i$ be the set of generators of type $i = 1,2,3$. We have $|\mc{G}_1| = 6, |\mc{G}_2|=5, |\mc{G}_3| = 5$. The binary variables $u_{g_i}, g_i\in\mc{G}_i, i = 1,2,3$, represent startup decisions, and the continuous variables $p_{g_i}, g_i\in\mc{G}_i, , i = 1,2,3$, represent production decisions. Scarf's example solves the following cost minimization problem:
\begin{subequations}\label{eq: scarf}
\begin{alignat}{4}
\min~~& \sum_{g_1\in\mc{G}_1} (53 u_{g_1} + 3p_{g_1}) + \sum_{g_2\in\mc{G}_2} (30u_{g_2} + 2p_{g_2}) + \sum_{g_3\in\mc{G}_3} 7 p_{g_3}\label{eq: scarf_obj}\\
\st~~& \sum_{g_1\in\mc{G}_1} p_{g_1} + \sum_{g_2\in\mc{G}_2} p_{g_2} + \sum_{g_3\in\mc{G}_3} p_{g_3} = d &&\label{eq: scarf_1} \\
& p_{g_3} \geq 2 u_{g_3} &&\forall g_3\in\mc{G}_3 \label{eq: scarf_2}\\
& p_{g_1} \leq 16 u_{g_1} &&\forall g_1\in\mc{G}_1 \label{eq: scarf_3}\\
&  p_{g_2} \leq 7 u_{g_2} &&\forall g_2\in\mc{G}_2 \label{eq: scarf_4}\\
& p_{g_3} \leq 6 u_{g_3} &&\forall g_3\in\mc{G}_3 \label{eq: scarf_5}\\
& p_{g}\geq 0, u_{g} \in\{0,1\} && \forall g\in\mc{G}_1, \mc{G}_2, \mc{G}_3\label{eq: scarf_6}
\end{alignat}
\end{subequations}
where the objective is to minimize the total cost. Constraint \eqref{eq: scarf_1} ensures the total production equals the demand. Constraints \eqref{eq: scarf_2} set the lower bound for the production of medium technology generators when they are on. Constraints \eqref{eq: scarf_3}-\eqref{eq: scarf_5} set the capacity of each generator. 

\begin{figure*}[h]
\centering
\begin{subfigure}[t]{0.4\textwidth}
\centering
\includegraphics[width = \textwidth]{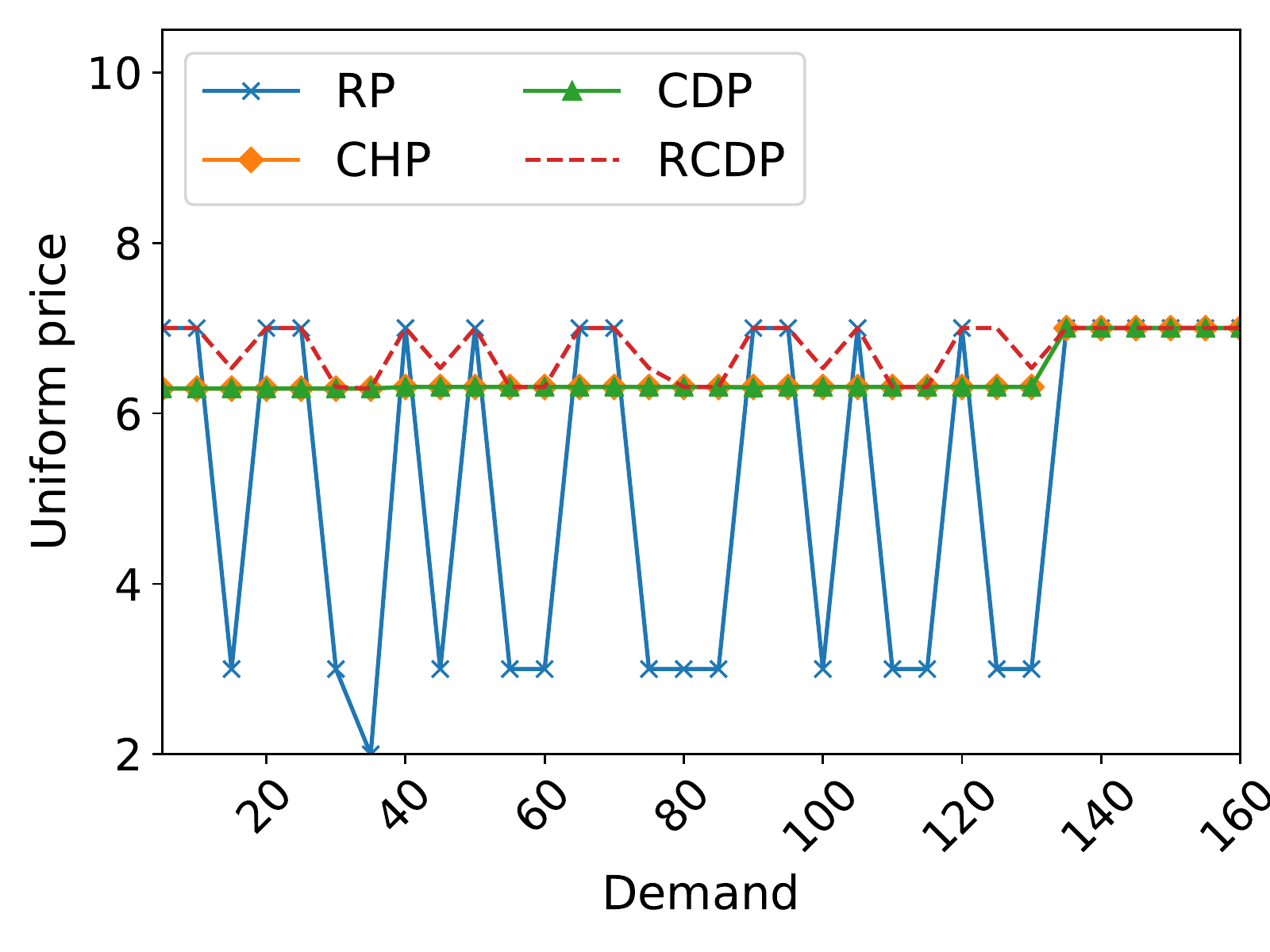}
\caption{}\label{fig: scarf_price}
\end{subfigure}
~
\begin{subfigure}[t]{0.4\textwidth}
\centering
\includegraphics[width = \textwidth]{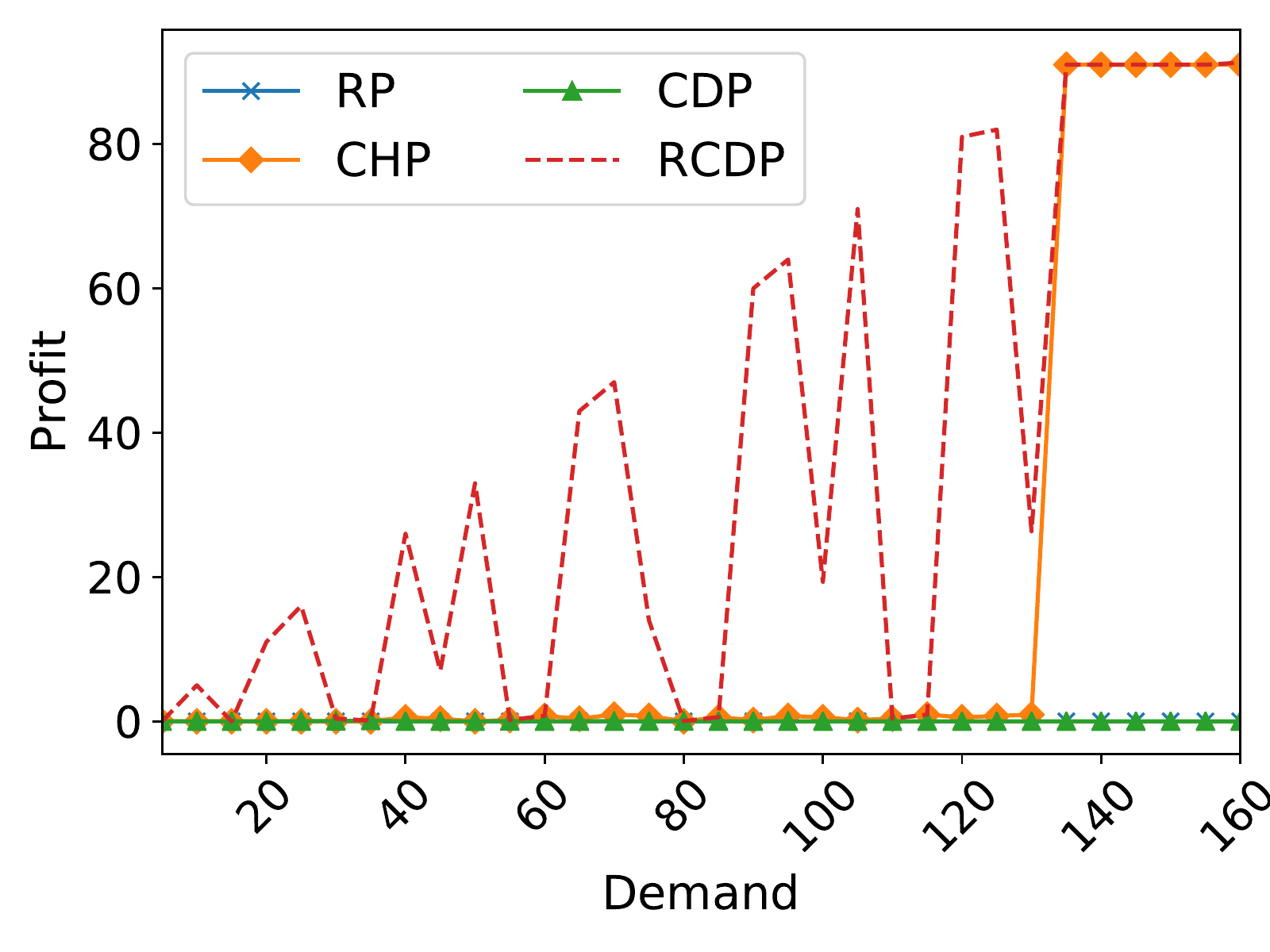}
\caption{}\label{fig: scarf_profit}
\end{subfigure}
~
\begin{subfigure}[t]{0.4\textwidth}
\centering
\includegraphics[width = \textwidth]{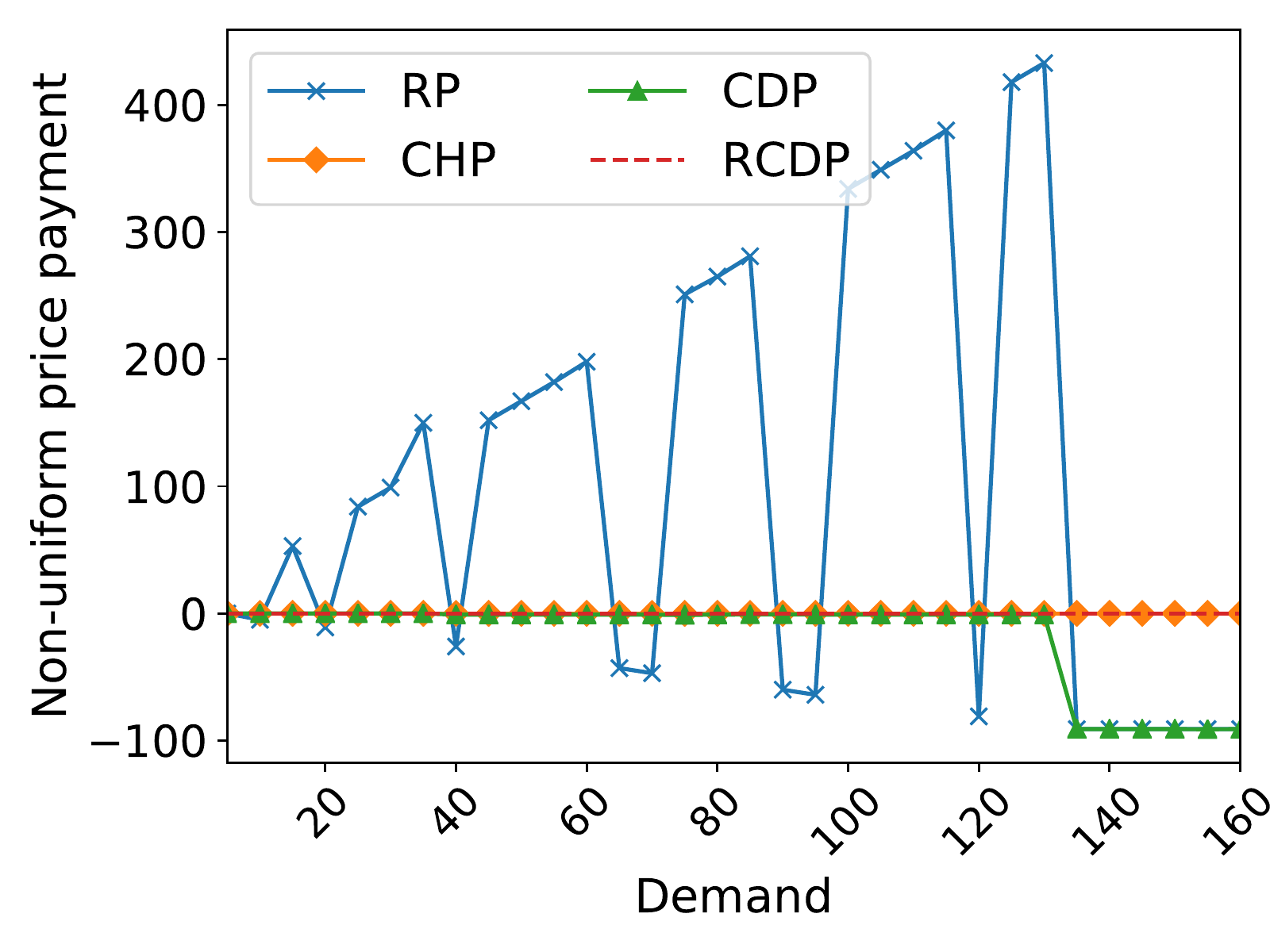}
\caption{}\label{fig: scarf_genDep}
\end{subfigure}
~
\begin{subfigure}[t]{0.4\textwidth}
\centering
\includegraphics[width = \textwidth]{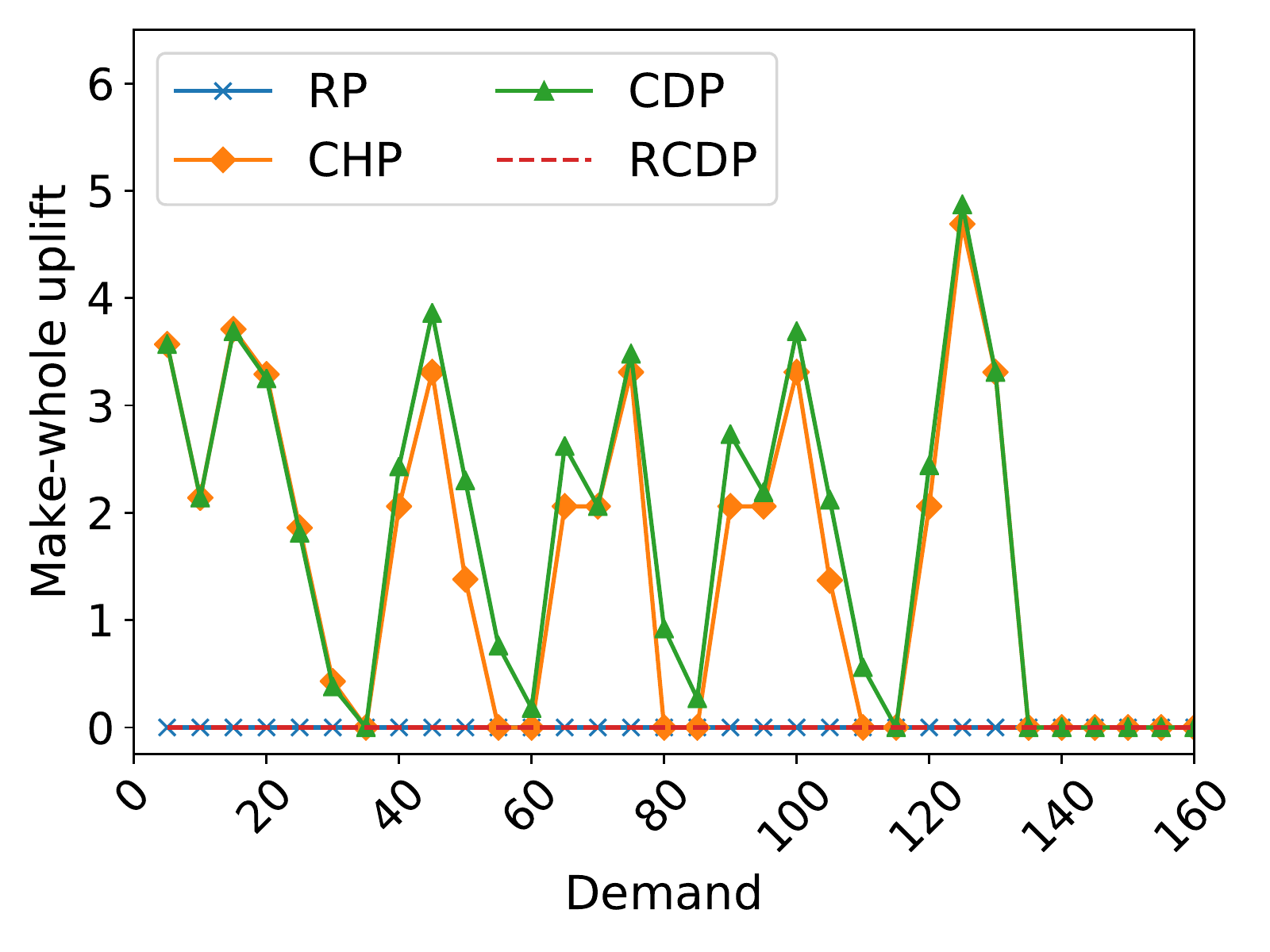}
\caption{}\label{fig: scarf_makewhole}
\end{subfigure}
\caption{\small Comparison of different pricing schemes for (a) uniform prices, (b) profits, (c) payments from generator-dependent prices, (d) make-whole uplift payments.}
\label{fig: scarf}
\end{figure*}

We experiment with various demand levels from 5 to 160, with a step length of 5. In Figure \ref{fig: scarf} we compare the following aspects of RP, CHP, CDP pricing and RCDP:

\begin{enumerate}
\item The uniform prices. Notice that for CDP and RCDP, both $\lambda_t^*$ and $\Lambda_t^*$ are uniform prices. In all of our experiments $\Lambda_t^*$ equals zero. Therefore, we report only $\lambda_t^*$ for those schemes.

\item Generator profits, as calculated by deducting costs from total revenue, which includes both price-based payments and uplift.

\item The payments from generator-dependent (non-uniform) prices $\phi^*_{jgt}$ and $\Phi^*_{jgt}$.

\item The make-whole uplift payments. This payment is made when the revenue from electricity prices is not enough to cover the costs. It is equal to the difference between the revenue and costs.
\end{enumerate}

Mosek solves all instances in under ten seconds. Because the cutting plane algorithm is significantly slower, we only use it when the demand level is less than 100. When the demand level is higher, we present results from the SDP approximation. 

Figure \ref{fig: scarf_price} shows that a small change in demand level can result in significant volatility in RP. This observation is consistent with results in the literature. Interestingly, CHP and CDP are equivalent for all demand levels. RCDP is higher than CDP for lower demand, and equals COP when the demand is high.

In Figure \ref{fig: scarf_profit}, we find that RP and CDP have zero profit for all instances. CHP generates near-zero profits at lower demand level and higher profits at higher demand levels. RCDP generates the highest profits among all pricing schemes, and match the profits of CHP at higher demand levels.

In Figure \ref{fig: scarf_genDep}, both CHP and RCDP have no generator-dependent pricesdue to the fact that both only use uniform prices corresponding to the demand constraint. CDP produces near-zero negative generator-dependent prices at low demand levels, and more negative prices at higher demand levels. RP produces volatile and large generator-dependent prices in many instances. As explained by \cite{o2005efficient}, the negative generator-dependent prices are used to discourage the entry of marginal plants when it is uneconomic to do so. In practice, utilities usually disregard such negative prices. 

Figure \ref{fig: scarf_makewhole} shows that RP requires zero make-whole payment, which is also consistent with the results in \cite{azizan2020optimal}. RCDP ensures revenue adequacy and thus also needs no make-whole payment. CHP requires make-whole payments, as expected, because Lagrangian duals of MIPs do not in general have strong duality. Interestingly, CDP also requires make-whole payments in many instances, as in those instances strong duality does not hold. 

The advantages of RCDP are that it does not rely on generator-dependent/make-whole uplift payments to ensure revenue adequacy, and its uniform prices are less volatile than RP.

\subsection{Pricing in Electricity Markets}\label{subsec: electricity_markets_results}
In this section we compare pricing schemes for unit commitment, as described in Section \ref{ch: app_pricing}. Note that in the experiments we use the formulation in Example \ref{case: simple_uc} for the $\uc$ problem. Our test instances are based on the adapted California ISO dataset of \cite{guoageneration}. The parameters for the generators are listed in Table \ref{tab:para_uc}.
\begin{table}[h]
\centering
\caption{Generator Parameters}
\label{tab:para_uc}
\begin{tabular}{lrrrr}
\hline
                  & \multicolumn{1}{c}{Gen. 1} & \multicolumn{1}{c}{Gen. 2} & \multicolumn{1}{c}{Gen. 3} & \multicolumn{1}{c}{Gen. 4} \\ \hline
$c_g^o$           & 25.0                       & 25.5                       & 44.68                      & 44.68                      \\
$c_g^s$           & 140.94                     & 140.94                     & 86.31                      & 86.31                      \\
$\underline{P}_g$ & 297                        & 238                        & 198                        & 198                        \\
$\overline{P}_g$  & 620                        & 496                        & 620                        & 620                        \\ \hline
\end{tabular}
\end{table}

We first consider 2 generators over 4 hours, which we refer to as Case 1. We set $d_t = [508, 644, 742, 776]$ and use Generators 1 and 2 in Table \ref{tab:para_uc}. We use the cutting plane algorithm of Section \ref{ch: cut_plane} to solve the CDP and RCDP COPs. Note that in this case, Mosek fails to solve the SDP approximation of the RCDP problem (returns ``UNKNOWN\_RESULT\_STATUS"), while the cutting plane algorithm solves it optimally in 282.85 seconds. When using the cutting plane algorithm, we set the bound for each element of the copositive matrix to 1000, and add complementary slackness constraint $\tr(\xstar \mathbf{x}^{*\top} \Omega) = 0$ when solving the CDP COP. 

We compare the following aspects of RP, CHP, CDP and RCDP and present the results in Table \ref{tab:uc_compare_1}:

\begin{enumerate}
\item Generator profits.

\item Profits without the uplift payments.

\item The payments from generator-dependent prices.

\item The make-whole uplift payments. 
\end{enumerate}

\begin{table}[h]
\centering
\caption{Comparison of Pricing Schemes for Case 1}
\label{tab:uc_compare_1}
\begin{tabular}{lrrrr}
\hline
                & \multicolumn{1}{c}{RP} & \multicolumn{1}{c}{CHP} & \multicolumn{1}{c}{CDP} & \multicolumn{1}{c}{RCDP} \\ \hline
Gen. 1 profit         & 0                      & 839.5                  & 292.0                  & 1393.6                  \\
Gen. 1 profit (pre-uplift)         & 0                      & 839.5                  & 292.0                  & 1393.6                  \\
Gen. 2 profit         & 0                      & 0                       & 0                       & 72.1                    \\
Gen. 2 profit (pre-uplift)         & 0                      & -95.3                       & -292.0                       & 72.1                    \\
Total profit     & 0                      & 839.5                  & 292.0                  & 1465.7                  \\
Total profit (pre-uplift)     & 0                      & 744.2                  & 0                  & 1465.7                  \\
Gen. dep. payment    & 498.0                 & 0                       & -201.5                 & 0                        \\
Make-whole uplift & 0                      & 95.3                   & 292.0                  & 0                        \\ \hline
\end{tabular}
\end{table}

RP results in zero profit for both generators, but relies on the revenue from generator-dependent prices for Generator 2. CHP uses a make-whole uplift payment to avoid a loss for Generator 2. Interestingly, CDP results in positive profit for Generator 1, and uses a make-whole uplift equal to Generator 1's profit to cover the loss of Generator 2. As expected, CDP is revenue neutral in aggregate, but not for individual generators. RCDP is the only pricing scheme that has no generator-dependent or make-whole uplift payments. It also ensures the generators with lower costs receive higher profits, which is desirable and to be expected in a market with only uniform prices.

In the second example, we include all 4 generators in Table \ref{tab:para_uc}, and 4 hours in the planning horizon. We refer to this as Case 2. The demand is $d_t = [1469, 1862, 2144, 2242]$. In Table \ref{tab:uc_compare_2}, we again compare the total profit, generator dependent payment, and make-whole uplift for each pricing scheme. In this case, Mosek is unable to solve the SDP approximation for the CDP or RCDP COPs, and the cutting plane algorithm does not terminate in a short time. We therefore use the result of the cutting plane algorithm after 3,000 iterations. Although the prices from CDP and RCDP are not optimal, we observe that the total revenue from uniform and availability prices, as well as the profit, decrease monotonically as more cuts are added (see e-companion Section \ref{ec: trend_graph}). This indicates that adding more cutting planes brings the prices closer to satisfying the theoretical properties of CDP and RCDP.

When using the cutting plane algorithm, we set the bounds on the absolute values of the elements in the copositive matrices to 5000. 
\begin{table}[h]
\centering
\caption{Comparison of Pricing Schemes for Case 2}
\label{tab:uc_compare_2}
\begin{tabular}{lrrrr}
\hline
                & \multicolumn{1}{c}{RP} & \multicolumn{1}{c}{CHP} & \multicolumn{1}{c}{CDP} & \multicolumn{1}{c}{RCDP} \\ \hline
Total profit     & 86859.5               & 86528.4                &  56560.1                       &129918.3                          \\
Total profit (pre-uplift)     & 86859.5               & 86445.0               &  56560.1                       &129918.3                          \\
Gen. dep. payment    & 86.3                  & 0                       &  -89354.5                       & 0                         \\
Make-whole uplift & 0                      & 83.4                   &   0                      & 0                         \\ \hline
\end{tabular}
\end{table}

The results for Case 2 are similar to Case 1. RCDP is the only pricing schemes that does not result in generator-dependent or make-whole uplift payments. It also results in larger profits than the other schemes. Although not solved to optimality, CDP still produces lower overall profit than RP and CHP. This is to be expected because the strong duality of COP guarantees zero pre-uplift profit. Given enough time to converge, CDP result in zero profit, excluding the make-whole uplift.

\subsection{Bimatrix Games}\label{subsec: bimatrix}
In this section, we provide a proof of concept for the copositive formulation of the KKT conditions, \eqref{eq: kkt_mip}, with bimatrix games. A bimatrix game has two players, 1 and 2, and each player has $n_i$ strategies to choose from. We use the binary vector $\x_i\in\{0,1\}^{n_i}$ to denote player $i$'s decision. The elements of $\x_i$ are binary because we only look for PNE. Player $i$ faces the following maximization problem:
\begin{subequations}\label{eq: bimatrix}
\begin{alignat}{4}
\max~~ &\x_1^\top R_i \x_2\label{eq: bimatrix_obj}\\
\st ~~& \sum_{j=1}^{n_i} x_{ij} = 1\label{eq: bimatrix_1}\\
& x_{ij} \in\{0,1\} &~~&\forall j=1,...,n_i,
\end{alignat}
\end{subequations}
where $R_i\in\mc{R}^{n_1\times n_2}$ is the payoff matrix. If Player 1 plays the $k_1^{\textrm{th}}$ strategy and Player 2 plays the $k_2^{\textrm{th}}$ strategy, then the payoff for Player $i$ is $R_{i,{k_1 k_2}}$. Constraint \eqref{eq: bimatrix_1} ensures that each player can only play one strategy.

Problem \eqref{eq: bimatrix} has all binary variables and only equality constraints. The upper bound of binary variables are implied by constraint \eqref{eq: bimatrix_1}. We therefore don't need to add upper bounds to the CPP reformulation. Because there are no quadratic terms between Player $i$'s decisions in the objective, $Q^{(i)}_{ii}$ in \eqref{eq: cpp_game_obj} of $\mc{P}_i^\cpp(\xni)$ is a zero matrix, which implies that $Q^{(i)}_{ii}\succeq 0$. Therefore, bimatrix games satisfy the requirements of the purely binary case in Section \ref{subch: compute_ne}. We can use the KKT conditions \eqref{eq: kkt_mip} to solve for the PNE of its CP counterpart, and then recover the PNE of the bimatrix game from the PNE of the CP game. 

We test the KKT conditions on bimatrix games with 2 to 5 strategies. For each number of strategies we test 5 instances. For all cases we use the cutting plane algorithm to solve KKT conditions \eqref{eq: kkt_mip}. In each initial master problem, we set an upper bound of 100 for the absolute values of the elements of each copositive matrix. Our computational results show that this bound is large enough to obtain copositive matrices upon termination.
We verify the correctness of our results by comparing with the results from the lrsnash bimatrix game solver \citep{avis2010enumeration}. {Note that our method is slower than specialized bimatrix game algorithms such as in \cite{avis2010enumeration}, but is applicable to more general classes of games.}

In Table \ref{tab: bimtrix} we {show the average time (seconds) and number of iterations used by the cutting plane algorithm in each case.} All test cases are solved in 7 seconds and under 90 iterations.

\begin{table}[h]
\centering
\caption{Performance of KKT Conditions for Bimatrix Games}
\label{tab: bimtrix}
\begin{tabular}{crr}
\hline
\#Strategies & \multicolumn{1}{c}{Time} & \multicolumn{1}{c}{Iterations} \\ \hline
2          & 1.5                           & 3                                 \\
3          & 1.5                           & 1                                   \\
4          & 2.5                           & 31                                  \\
5          & 4.8                          & 62                                  \\ \hline
\end{tabular}
\end{table}
\section{Conclusion}
\label{sec:conc}
MBQPs can be equivalently written as CPPs, which are NP-hard but convex. Given an MBQP, we straightforwardly derive its dual COP. Due to convexity, if a constraint qualification is satisfied, the CPP and COP have strong duality. This provides a new and general notion of duality for discrete optimization problems.

We apply this perspective in two ways. We first design a new pricing mechanism for nonconvex electricity markets, which has several useful theoretical properties. One direction of future study is the design of economic mechanisms for other nonconvex markets, e.g., surge pricing in transportation. Second, we reformulate MBQ games as CP games and provide conditions for the existence and uniqueness of equilibria. We also use the KKT conditions to solve for the equilibria. To enable implementation, we design a new cutting plane algorithm for COPs, which we use to solve our numerical examples.

There are several promising avenues of future work. It would be useful to identify classes of MBQP that, when reformulated as CPPs, always have strong duality. By reformulating other classes of discrete games in terms of CPPs, existing results for convex games could potentially identify new conditions for existence and uniqueness. We also intend to extend our results to mixed-strategy Nash equilibria. Finally, it may be possible to improve the cutting plane algorithm by strengthening the master problem, and by deriving conditions under which the algorithm is guaranteed to terminate.


%% file: ecompanion.tex
\section{SDP Approximation for COP}\label{ec: sdp_approx}
Let $\mc{S}_n^+$ be the $n$-dimensional positive semidefinite (PSD) cone and $\mc{N}_n$ be the cone of $n$-dimensional entrywise nonnegative matrices, then we have $\mc{C}_n^*\subseteq \mc{S}_n^+ \cap \mc{N}_n$ and $\mc{S}_n^+ + \mc{N}_n\subseteq \mc{C}_n$ \citep{dur2010copositive}. Taking advantage of these relationships, we can relax and approximately solve CPP as an optimization problem over the intersection of $\mc{S}_n^+$ and $\mc{N}_n$ cones, and restrict the COP problem as an optimization over the Minkowski sum of $\mc{S}_n^+$ and $\mc{N}_n$ cones. 

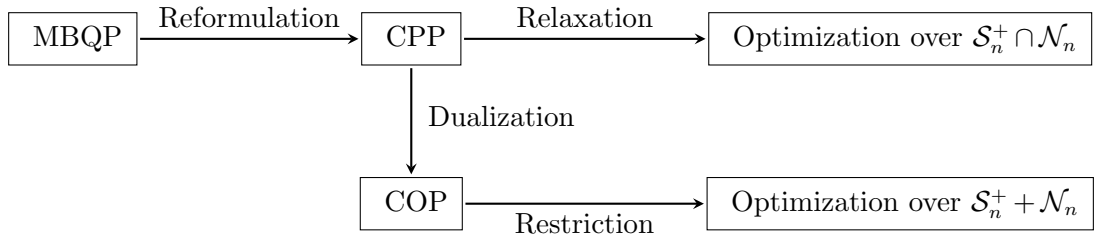
\begin{figure}[h]
\centering
\begin{tikzpicture}[arrow/.style={thick,->,shorten >=1pt,shorten <=1pt,>=stealth},]
	\tikzstyle{block} = [rectangle, draw, text centered, minimum height=1.8em, minimum width=0.5em]
	\node[right of=A,xshift = 3.5cm,block,midway,above,yshift = 0.5em] (B) {
    		 \renewcommand\arraystretch{0.5} $\begin{array}{c} \text{CPP} \end{array}$};
	\node[below of=B,xshift = 4.5cm,block,midway,above,yshift = -1cm] (C) {
    		 \renewcommand\arraystretch{0.5} $\begin{array}{c} \text{COP} \end{array}$};
        \node[block,midway,above,yshift = 0.5em] (A) {
    		 \renewcommand\arraystretch{0.5} $\begin{array}{c} \text{MBQP}  \end{array}$};
	\node[right of=B,xshift = 10cm,block,midway,above,yshift = 0.5em] (rel_cpp) {
    		 \renewcommand\arraystretch{0.5} $\begin{array}{c} \text{Optimization over $\mc{S}_n^+ \cap \mc{N}_n$} \end{array}$};
	\node[right of=C,xshift = 10cm,block,midway,above,yshift = -2cm] (rst_cop) {
    		 \renewcommand\arraystretch{0.5} $\begin{array}{c} \text{Optimization over $\mc{S}_n^+ + \mc{N}_n$} \end{array}$};
		 
	\draw[arrow] (A.east) to node[sloped,swap,anchor=center, above][yshift=0.5pt] {Reformulation} (B.west); 
	\draw[arrow] (B.south) to node[swap,anchor=center, above][yshift=-2mm, xshift = 1.2cm] {Dualization} (C.north); 
	\draw[arrow] (B.east) to node[sloped,swap,anchor=center, above][yshift=0.5pt] {Relaxation} (rel_cpp.west); 
	\draw[arrow] (C.east) to node[sloped,swap,anchor=center, below][yshift=0.5pt] {Restriction} (rst_cop.west); 
\end{tikzpicture}
\caption{Relationships between MBQP, CPP, COP and their approximations}\label{fig: cop_relation}
\end{figure}

In our work, we reformulate MBQP models to CPPs using the method by \cite{burer2009copositive}, to utilize the dual COP problems. We summarize the relationship between all those mathematical programming problems in Figure \ref{fig: cop_relation}.

We need to solve COPs in this work. One method often used in the literature, as mentioned above, is to use the relationship $\mc{S}_n^+ + \mc{N}_n\subseteq \mc{C}_n$ for approximation. More specifically, we replace the conic constraint $\Omega\in\mc{C}_n$ with the following restriction: 
\begin{align*}
\Omega - N &\in\mc{S}_n^+\\
N &\geq 0,
\end{align*}
which can be solved with SDP solvers such as Mosek and SeDuMi.

Another method to obtain a solution of a COP is to solve its dual CPP problem using a commercial solver, then query the duals of CPP constraints via the solver. However, there is not any solver that directly solves CPPs, so we instead solve an SDP relaxation of the CPP problem, then query the duals of the SDP relaxation. More specifically, we relax the conic constraint $X\in\mc{C}_n^*$ to the following constraints:
\begin{align*}
X &\in\mc{S}_n^+\\
X &\geq 0,
\end{align*}
which can then be solved with SDP solvers.
\section{Symmetry of CPP}\label{ec: sym_cpp}
When taking of dual of conic programs, it is often assumed in the literature that the coefficient matrices are symmetric. We find that this assumption may lead to a difference in the form of the dual problem, if it is also assumed that the conic matrix is in $\mc{S}_n$, the set of symmetric matrices. We include a discussion about this issue in this section, as we believe it is an important detail to notice when deriving dual problems in copositive programming, and it is not something we find discussed in the literature.

We study the following CPP: 
\begin{subequations}\label{eq: ec_cpp}
\begin{alignat}{4}
\mc{P}^\cpp_{\mc{C}}: \min ~~&\tr(C^\top X)~~ \label{eq: ec_cpp_obj}\\
\st ~~&\tr(A_i^\top X) = b_i&~~&\forall i = 1, ..., m &~~& (\lambda_i) \label{eq: ec_cpp_1}\\
&X\in\mc{C}_n^*&& &~~&(\Omega). \label{eq: ec_cpp_2}
\end{alignat}
\end{subequations}

If $C\in\mc{S}_n$ and $A_i\in\mc{S}_n, i=1,...,m$, then the dual of $\mc{P}^\cpp_{\mc{C}}$ is the following:
\begin{subequations}\label{eq: ec_cop}
\begin{alignat}{4}
\mc{P}^\cop_{\mc{C}}: \min ~~&\lambda^\top b ~~ \label{eq: ec_cop_obj}\\
\st ~~&C-\sum_{i=1}^m \lambda_i A_i -\Omega = 0 \label{eq: ec_cop_1}\\
&\Omega\in\mc{C}_n.&& &~~&& \label{eq: ec_cop_2}
\end{alignat}
\end{subequations}

However, if some of the parameter matrices are not symmetric, $\mc{P}^\cop_{\mc{C}}$ is not necessarily the correct dual formulation. To understand this, first notice that although $\mc{C}_n^*$ contains only symmetric matrices, its dual cone in $\mb{R}^{n\times n}$ may contain non-symmetric matrices\footnotetext{A similar discussion for PSD cones can be found at \url{https://github.com/cvxgrp/scs/issues/31}.}. For example, the non-symmetric matrix $Y = \begin{bmatrix}10 & 1\\-1 & 10\end{bmatrix}$ is in the dual cone of $C_n^*$. This is because for any $X\in\mc{C}_2^*$, 
$$\tr(XY) = \tr\left(X \begin{bmatrix}10 & 0\\0 & 10\end{bmatrix}\right) + \tr\left(X \begin{bmatrix}0 & 1\\-1 & 0\end{bmatrix}\right) = 10\tr(X) + 0 \geq 0.$$ 
where the last inequality follows from the fact that when $n\leq 4$, $\mc{C}_n^* = \mc{S}_n^+\cap\mc{N}_n$ \citep{dur2010copositive}\footnotemark.

To derive the dual problem of $\mc{P}^\cpp_{\mc{C}}$ in the general case, we define general completely positive cone and copositive cone that contain non-symmetric matrices, by using their symmetric parts:
$$\mc{GC}_n^*:=\left\{X\in\mb{R}^{n\times n}\;\Bigg|\; \frac{X + X^\top}{2}=\sum_{k\in\mc{K}} z^k (z^k)^\top \text{for some finite~} \{z^k\}_{k\in \mc{K}}\subset \mb{R}^{n}_+\sm \{\mathbf{0}\} \right\},$$
and 
$$\mc{GC}_n := \left\{X\in\mb{R}^{n\times n}\;\Bigg|\; \mathbf{y}^\top \frac{X + X^\top}{2} \mathbf{y} \geq0\textrm{ for all } \mathbf{y}\in\mathbb{R}^n_+\right\}.$$
It can be proved that $\mc{GC}_n^*$ and $\mc{GC}_n$ are dual cones of each other. 

Then $\mc{P}^\cpp_{\mc{C}}$ is equivalent to the following CPP:
\begin{subequations}\label{eq: gc_cpp2}
\begin{alignat}{4}
\mc{P}^{\cpp}_{\mc{GC}}: \min ~~&\tr(C^\top X)~~ \label{eq: gc_cpp2_obj}\\
\st ~~&\tr(A_i^\top X) = b_i&~~&\forall i = 1, ..., m && (\lambda_i) \label{eq: gc_cpp2_1}\\
&X\in\mc{GC}_n^* && &&(\Omega) \label{eq: gc_cpp2_2}\\
&X_{ij}= X_{ji} &&\forall i=1,...,n-1, \forall j = i+1, ...,n ~~&&(\pi_{ij}). \label{eq: gc_cpp2_3}
\end{alignat}
\end{subequations}

Thus, the dual of both $\mc{P}^\cpp_{\mc{C}}$ and $\mc{P}^{\cpp}_{\mc{GC}}$ is:
\begin{subequations}\label{eq: gc_cop2}
\begin{alignat}{4}
\mc{P}^{\cop}_{\mc{GC}}: \min ~~&\lambda^\top b ~~ \label{eq: gc_cop2_obj}\\
\st ~~&C-\sum_{i=1}^m \lambda_i A_i + \sum_{i=1}^{n-1}\sum_{j=i+1}^n \pi_{ij} (e_i e_j^\top - e_j e_i^\top) -\Omega = 0 \label{eq: gc_cop2_1}\\
&\Omega\in\mc{GC}_n.&& &~~&& \label{eq: gc_cop2_2}
\end{alignat}
\end{subequations}
We have the following proposition:
\begin{proposition}
If matrices $C$ and $A_i, i=1,...,m$ are symmetric, then $\mc{P}^{\cop}_{\mc{GC}}$ and $\mc{P}^\cop_{\mc{C}}$ are equivalent.
\end{proposition}
\proof{Proof.}
Let $\hat{X}$ be an optimal solution of $\mc{P}_{\mc{C}}^\cpp$ and thus also of $\mc{P}_{\mc{GC}}^\cpp$. Since $C$ and $A_i, i=1,...,m$ are symmetric, $\tilde{X} := (\hat{X}+\hat{X})/2$ is also optimal to the optimization problem \eqref{eq: gc_cpp2_obj}-\eqref{eq: gc_cpp2_2} ($\mc{P}_{\mc{GC}}^\cpp$ without the symmetry constraints). 
Thus at $\tilde{X}$, constraints \eqref{eq: gc_cpp2_3} are redundant, which means there exists an optimal dual solution with all $\pi_{ij} = 0$.

This means in $\mc{P}^{\cop}_{\mc{GC}}$ we can get rid of the term $\sum_{i=1}^{n-1}\sum_{j=i+1}^n \pi_{ij} (e_i e_j^\top - e_j e_i^\top)$, while preserving at least one optimal solution. This makes $\mc{P}^{\cop}_{\mc{GC}}$ the same as $\mc{P}^\cop_{\mc{C}}$ except for the conic constraint being non-symmetric. However, the symmetry of $\Omega$ in $\mc{P}^{\cop}_{\mc{GC}}$ is ensured by the symmetry of the matrices in constraint \eqref{eq: gc_cop2_1}. Thus, $\mc{P}^{\cop}_{\mc{GC}}$ is equivalent to $\mc{P}^\cop_{\mc{C}}$. 
\Halmos\endproof

From the proof we can conclude that if some $A_i$ is not symmetric, then $\mc{P}^\cop_{\mc{C}}$ is generally not the dual of $\mc{P}^{\cpp}_{\mc{C}}$.
\section{Proofs}\label{ec: proofs}
\begin{repeattheorem}[Proposition \ref{th: cop_cut}.]
If the optimal value of $\SP(\bar{\Omega})$ is nonzero, then \eqref{eq: cop_cut} cuts off $\bar{\Omega}$.
\end{repeattheorem}
\proof{Proof of Proposition \ref{th: cop_cut}.}
Because of \eqref{eq: separation_c} and the fact that $q=1$, there is at least one element in $u$ that is nonzero, i.e. $\beta = \{i | u_i = 1, i = 1, ...,n_c\} \neq \emptyset$. Denote $\Omega_{\beta\beta}$ as the submatrix of $\Omega$ that consists of rows and columns with indices in $\beta$. Similarly we define the subvector $\mathbf{z}_{\beta}$. From the fact that the optimal objective $\bar{w}>0$, we know from constraint \eqref{eq: separation_b} that $\bar{\Omega}_{\beta\beta}\bar{\mathbf{z}}_\beta < \mathbf{0}$ and thus $\bar{\mathbf{z}}_\beta \neq \mathbf{0}$. Therefore, $\bar{\mathbf{z}}_\beta^\top\bar{\Omega}_{\beta\beta}\bar{\mathbf{z}}_\beta < 0$. Also, let $\alpha = \{1,...,n_c\} \sm \beta$, then $\bar{\mathbf{z}}_\alpha = \mathbf{0}$ due to \eqref{eq: linkingZandU}, which means $\bar{\mathbf{z}}^\top\bar{\Omega}\bar{\mathbf{z}} = \bar{\mathbf{z}}_\beta^\top\bar{\Omega}_{\beta\beta}\bar{\mathbf{z}}_\beta < 0$. Thus, $\bar{\Omega}$ violates the cut \eqref{eq: cop_cut}. 
\Halmos\endproof

\begin{repeattheorem}[Theorem \ref{th: rev_adequacy}.]
If strong duality holds for $\uc^\cpp$, then CDP balances the revenue and the aggregate cost of the generators. It also balances the revenue and payment of the SO.
\end{repeattheorem}
\proof{Proof of Theorem \ref{th: rev_adequacy}.}
Fix all primal and dual variables at the optimal values. Multiplying constraints \eqref{eq: uc1_cpp1_1} - \eqref{eq: uc1_cpp1_8} by their corresponding dual variables yields
\begin{subequations}\label{eq: uc1_balance}
\addtocounter{equation}{1}
\begin{alignat}{4}
& \sum_{g\in\mc{G}}\lambda_t^* p^*_{gt} = \lambda_t^* d_t && \forall t\in\mc{T} \label{eq: uc1_balance_1}\\
&\phi_{jgt}^* \aphitop_{jgt} \x^* = \phi_{jgt}^* b_{jgt} &&\forall  j= 1,...,m, g\in\mc{G}, t\in\mc{T} \label{eq: uc1_balance_ub}\\
&\sum_{g\in\mc{G}} \Lambda^*_t \Xpstar + 2\sum_{g_1< g_2, g_1,g_2 \in\mc{G}}\Lambda^*_t \Xppstar = \Lambda^*_t d_t^2&~~~&\forall t\in\mc{T}\\
&\Phi_{jgt}^* \tr(\aphi_{jgt} \aphitop_{jgt} X^*) = \Phi_{jgt}^* b_{jgt}^2  &~&\forall j= 1,...,m, g\in\mc{G}, t\in\mc{T}.
\end{alignat}
\end{subequations}

Summing all the constraints in \eqref{eq: uc1_balance}, we obtain 
\begin{align}
\label{eq: payment_dual_obj}
&  \sum_{t\in\mc{T}}\sum_{g\in\mc{G}}\left(\lambda_t^* p^*_{gt} + \Lambda^*_t \Xpstar+ \sum_{j=1}^m\left(\phi^{*}_{jgt}\aphi_{jgt} \x^* + \Phi^{*}_{jgt}\tr(\aphi_{jgt} \aphitop_{jgt} X^*)\right) \right) \nonumber\\&+ 2\sum_{g_1< g_2, g_1,g_2 \in\mc{G}}\Lambda^*_t \Xppstar \nonumber\\
 =& \sum_{t\in\mc{T}}\left(d_t \lambda_t^*  + d_t^2 \Lambda_t^*  +  \sum_{j=1}^m \sum_{g\in\mc{G}}\left(b_{jgt}\phi_{jgt}^* + b_{jgt}^2 \Phi_{jgt}^* \right)\right)
\end{align}
where the left-hand side and the right-hand side are respectively equal to $\Pi^{\text G} := \sum_{g\in\mc{G}}\sum_{t\in\mc{T}}\pi^{\text G}_{gt}$ and $\Pi^{\text L} := \sum_{t\in\mc{T}}\pi^{\text L}_t$ (using the definitions of $\pi^{\text G}_{gt}$ and $\pi^{\text L}_t$ in Definition \ref{def:CDP}).
This shows that the total payment $\Pi^{\text G}$ is equal to the total revenue $\Pi^{\text L}$, thus CDP is revenue neutral for the SO.

Because we have assumed that strong duality holds for $\uc^\cpp$, and considering the optimal objective of $\uc^\cop$ equals $\pi^{\text L}$, we have:
\[
\pi^{\text G} = \sum_{g\in\mc{G}}\sum_{t\in\mc{T}} \left(c^p_g p^*_{gt} + c^u_g u^*_{gt}\right).
\]
Therefore, for the generators the total payment from the SO equals the total cost.
\Halmos\endproof
\begin{repeattheorem}[Theorem \ref{th: ind_rational}.]
Assume CPP \eqref{eq: uc1_cpp1} satisfies some sufficient condition for convex programming strong duality. Additionally, let $(\x^*, X^*)$ be an optimal solution for the problem \eqref{eq: uc1_cpp1}, assume all items in $X^*$ in the forms of $\Xvwstar,\;\forall g_1 \neq g_2, g_1,g_2\in\mc{G}, t\in\mc{T}$ are equal to zero. Then the market mechanism with CDP is also individually rational.
\end{repeattheorem}
\proof{Proof of Theorem \ref{th: ind_rational}.}
In the CPP reformulation \eqref{eq: uc1_cpp1}, we dualize demand constraints \eqref{eq: uc1_cpp1_1} and lifted demand constraints \eqref{eq: uc1_cpp1_7} with their respective optimal dual prices $\lambda_t^*$ and $\Lambda_t^*$, we obtain the following Lagrangian relaxation problem:
\begin{subequations}\label{eq: uc1_lagrangian}
\begin{alignat}{4}
&\min &&\sum_{g\in\mc{G}}\sum_{t\in\mc{T}} \left(c^o_g p_{gt} + c^s_g u_{gt}\right) + \sum_{t\in\mc{T}}  \lambda_t^* \left(d_t - \sum_{g\in\mc{G}}p_{gt}\right)\nonumber \\
& && +\sum_{t\in\mc{T}} \Lambda^*_t \left(d_t^2 - \sum_{g\in\mc{G}}\left(\Xp + \sum_{g'\in\mc{G}, g'\neq g} \Xpp\right)\right)\label{eq: uc1_lagrangian_1}\\
&\st && \eqref{eq: uc1_cpp1_2}, \eqref{eq: uc1_cpp1_8} - \eqref{eq: uc1_cpp1_14}.
\end{alignat} 
\end{subequations}

Since we have strong duality, and because the Lagrangian multipliers in the Lagrangian relaxation problem \eqref{eq: uc1_lagrangian} are fixed to their optimal values, an optimal solution $(\x^*, X^*)$ for CPP \eqref{eq: uc1_cpp1} is also optimal for its Lagrangian relaxation \eqref{eq: uc1_lagrangian}.

Because the terms in the form of $\Xvw$ are zeros, the term $\Xpp$ in the objective can be eliminated. Also, because of this condition, by Lemma 1 of \cite{drew1998completely}, $X$ is completely positive if and only if $X_g, \forall g\in\mc{G}$ are completely positive. Then straightforwardly we obtain the equivalence between $Y_g\in\mc{C}_{n_g+1}^*, \forall g\in\mc{G}$ and $Y\in\mc{C}_{n+1}^*$.

Now, the Lagrangian relaxation \eqref{eq: uc1_lagrangian} can be separated into individual optimization problems, one per generator. Converting the minimization problem into maximization, and ignoring the constant terms $\lambda^*_t d_t$ and $\Lambda^*_t d_t^2$ in the objective, solving problem \eqref{eq: uc1_lagrangian} is equivalent to solving the following maximization problem for all $g\in\mc{G}$:
\begin{subequations}\label{eq: uc1_individual_lagrangian}
\begin{alignat}{4}
\max~~ &\sum_{t\in\mc{T}} \left(\lambda^*_t p_{gt} + \Lambda^*_t \Xp - c^o_g p_{gt} - c^s_g u_{gt}\right)\label{eq: uc1_individual_lagrangian_1} \\
\st~~ &\aphitop_{jgt} \x = b_{jgt} \\
&\tr(\aphi_{jgt} \aphitop_{jgt} X) = b_{jgt}^2\\
&z_{gt} = Z_{gt} \\
& Y_g \in\mc{C}_{n_g+1}^*.
\end{alignat} 
\end{subequations}

Therefore, if a solution $(\x^*, X^*)$ is optimal for CPP \eqref{eq: uc1_cpp1}, then its component corresponding to $g$, $[\xstarg, X^{*}_g]$, also solves \eqref{eq: uc1_individual_lagrangian} optimally. 

Notice that in the profit-maximization problem \eqref{eq: uc1_individual} we assumed that at optimality $\Xpstar = p^{*2}_{gt}$ and $\Xppstar = p^*_{gt} p^*_{g' t}$, in addition to the assumption of $\Xvwstar = 0$, we can eliminate the cross-term payment from the the objective of \eqref{eq: uc1_individual}. Now if we reformulate \eqref{eq: uc1_individual} into a CPP, its objective is the same as \eqref{eq: uc1_individual_lagrangian_1} (except for the constant terms $\phi^*_{jgt}b_{jgt}$ and $\Phi^*_{jgt} b^2_{jgt}$). Moreover, problem \eqref{eq: uc1_individual_lagrangian} is exactly the CPP equivalence of the profit-maximization problem \eqref{eq: uc1_individual}, which means $\xstarg$ is optimal for \eqref{eq: uc1_individual}. Therefore, if we charge electricity prices $\lambda^*_t$ and $\Lambda^*_t$, individual generators will not have the incentive to deviate from the optimal dispatch levels set by the system operator, as those dispatch levels also maximize their own profits.
\Halmos\endproof
\begin{repeattheorem}[Theorem \ref{th: games_equiv}.]
Let $\mathbf{x}^*$ be a PNE of $\mc{G}^\bqp$. Then $(\xstari, \xstari \xstaritop)_{i\in\mc{I}}$ is a PNE of $\mc{G}^\cpp$. 

Conversely, if $(\xstari, \Xstari)_{i\in\mc{I}}$ is a PNE of $\mc{G}^\cpp$, and if $(\Qi_{ii}\succeq 0)_{i\in\mc{I}}$, $(\xstarik\in\{0,1\},\;\forall k\in\mc{B}_i)_{i\in\mc{I}}$, then $\mathbf{x}^*$ is a PNE of $\mc{G}^\bqp$. 
\end{repeattheorem}
\proof{Proof of Theorem \ref{th: games_equiv}.}
($\Rightarrow$) We first prove the conversion from $\mc{G}^\bqp$ to $\mc{G}^\cpp$. Since $\xstar$ is a PNE of $\mc{G}^\bqp$, it satisfies $f_i(\mathbf{x}^*_i, \mathbf{x}^*_{-i}) = \text{opt}(\mc{P}^\bqp_i(\xstarni)),\;\forall i \in\mc{I}$. On the other hand, $\Pith^\bqp(\xstarni)$ can be reformulated to 
$\Pith^\cpp(\xstarni)$. From Remark \ref{th: qip_sol}, we know that $(\xstari, \xstari \xstaritop)$ is an optimal solution of $\mc{P}^{\text{CPP}}_i(\xstarni)$, which means $(\xstari, \xstari \xstaritop)_{i\in\mc{I}}$ is a PNE for $\mc{G}^{\text{CPP}}$.

($\Leftarrow$) For the opposite direction, assume $(\xstari, X^*_i)_{i\in\mc{I}}$ is a PNE for $\mc{G}^\cpp$, which means $(\xstari, X^*_i)$ is an optimal solution for $\Pith^\cpp(\xstarni)$, and $\xstari$ is in the convex hull of optimal solutions for $\Pith^\bqp(\xstarni)$.

Now assume $(\Qi_{ii}\succeq 0)_{i\in\mc{I}}$, $(\xstarik\in\{0,1\},\;\forall k\in\mc{B}_i)_{i\in\mc{I}}$. 
We first show that $\xstari$ is feasible to $\Pith^\bqp(\xstarni)$. Note that  $\xstari$ satisfies \eqref{eq: ip_game_1} because it satisfies \eqref{eq: cpp_game_1}. $\xith\in\mb{R}^l_+$ is guaranteed because $\Xstari \in\mc{C}_{n_i}^* \Rightarrow \Xstari \geq 0  \Rightarrow X^*_i \alphai \geq 0  \Rightarrow \xstari \geq 0$, where the last step follows from the constraint \eqref{eq: cpp_game_4}. Additionally, we assumed $\xstarik\in\{0,1\},\;\forall k\in\mc{B}_i$, which means \eqref{eq: ip_game_3} is satisfied. Since $\xstari$ is feasible to $\Pith^\bqp(\xstarni)$ and $\Qi_{ii}\succeq 0$, by Remark \ref{th: feas_bqp_x} we know that $\xstari$ is an optimal solution for $\Pith^\bqp(\xstarni)$, where $\xstarni$ are also ensured to be pure strategies. Thus, $\xstar$ is a PNE for $\mc{G}^\bqp$. 
%
%
\Halmos\endproof

\begin{repeattheorem}[Proposition \ref{th: exist_qip}.]
{\bf (Existence of PNE)} $\mc{G}^\bqp$ has at least one PNE if
\BI
\I[(i)] Property \eqref{a: alpha} is satisfied for $\Pith^{\bqp}(\xni)$,

\I[(ii)] $\Qi_{ii}\succeq 0$, and

\I[(iii)] there exists a PNE $(\xstari, \Xstari)_{i\in\mc{I}}$ for the corresponding $\mc{G}^\cpp$ that satisfies $(\xstarik\in\{0,1\},\;\forall k\in\mc{B}_i)_{i\in\mc{I}}$.
\EI
\end{repeattheorem}
\proof{Proof of Proposition \ref{th: exist_qip}.}
Condition (i) ensures that the optimization problems of $\mc{G}^\bqp$ can be reformulated as CPPs: For each player $i$, because property \eqref{a: alpha} is true, we can reformulate $\Pith^{\bqp}(\xni)$ as the corresponding CPP problem $\Pith^{\cpp}(\xni)$. As we discussed in the text, the corresponding CP game $\mc{G}^\cpp$ has at least one PNE.

Let $(\xstari, \Xstari)_{i\in\mc{I}}$ be an PNE that satisfies condition (iii), along with the condition (ii), all requirements for the existence of PNE for $\mc{G}^\bqp$ in Theorem \ref{th: games_equiv} are met and thus $\xstar$ is a PNE for $\mc{G}^\bqp$. Therefore, the PNE of $\mc{G}^\bqp$ exists with conditions (i)-(iii).
\Halmos\endproof

\section{Individual Revenue Adequacy with Both Uniform and Availability Prices}\label{ec: ind_rev_avail_price}
If we include both uniform and availability prices in the revenue adequacy constraints, we have the following pricing problem to solve:
\begin{subequations}\label{eq: uc1_cop_profitable2}
\begin{alignat}{4}
\max~~&\sum_{t\in\mc{T}}\left( d_t\lambda_t +  d_t^2\Lambda_t+ \sum_{j=1}^m \sum_{g\in\mc{G}}(b_{jgt}\phi_{jgt}  +  b_{jgt}^2\Phi_{jgt})\right)\label{eq: uc1_cop_profitable2_obj}\\ 
\st ~~&\sum_{t\in\mc{T}} \left(p^*_{gt}\lambda_t + p^{*2}_{gt} \Lambda_t + \sum_{g'\in\mc{G}\sm\{g\}} p_{gt}^*p_{g't}^*\Lambda_t + \sum_{j=1}^m\left(\aphitop_j \xstar \phi_{jgt} + \tr(\aphi_j \aphitop_j X^*) \Phi_{jgt} \right) \right)\hspace{-1cm}\nonumber\\& \geq \sum_{t\in\mc{T}}\left(c^o_g p^*_{gt} + c^s_g u^*_{jgt}\right) &~~~&\forall g\in\mc{G} \label{eq: uc1_cop_profitable2_1}\\
&(\pmb \lambda, \pmb \phi, \pmb \Lambda, \pmb \Phi, \pmb \delta, \Omega)\in\mc{F}^{\text{COP}}.\hspace{-1cm}&&\label{eq: uc1_cop_profitable2_constr}
\end{alignat}
\end{subequations}
Again, we use $(p^{*}_{gt})^2$ in place of $\Xpstar$ and $p^*_{gt} p^*_{g' t}$ in place of $\Xppstar$.

If \eqref{eq: uc1_cop_profitable2} is feasible, then prices from \eqref{eq: uc1_cop_profitable2} should satisfy revenue neutrality for generators. This is because if we sum up left-hand sides of \eqref{eq: uc1_cop_profitable2_1} over $g$, then the value equals the objective \eqref{eq: uc1_cop_profitable2_obj} (proved by \eqref{eq: payment_dual_obj}), which represents the total revenue of generators. One the other hand, \eqref{eq: uc1_cop_profitable2} can be viewed as imposing extra constraints on the original CDP problem \eqref{eq: uc1_cop}, whose objective value, according to weak duality, is no more than the total costs $\sum_{g\in\mc{G}}\sum_{t\in\mc{T}} \left(c^o_g p^*_{gt} + c^s_g u^*_{gt}\right)$. But \eqref{eq: uc1_cop_profitable2} also restricts the total revenue of generators to be no less than the total costs. Therefore, a feasible solution of \eqref{eq: uc1_cop_profitable2} should ensure revenue neutrality for generators.

In addition, we actually have revenue neutrality for every generator, so each generator is paid exactly its cost. To understand why this result is true, assume towards contradiction if any generator has a strictly positive profit, then because of revenue neutrality of the whole system, some other generator must have a strictly negative profit, which violates \eqref{eq: uc1_cop_profitable2_1}. 

In comparison, we do not have revenue neutrality for individual generators in RCDP, i.e., in that case generators could have strictly positive profits.

Similar to the case of RCDP, with this pricing scheme, it can be proved that \eqref{eq: uc1_cop_profitable2} is guaranteed to be feasible if its dual problem has an interior.
\section{COP Algorithms Comparison for Maximum Clique Problem}\label{ec: max_clique} 
To compare our COP cutting plane algorithm with commonly used approaches for COP in literature, we present the results from solving the maximum clique problem with those algorithms. The maximum clique problem tries to find the maximum clique number on a graph $\mc{G} = (\mc{N}, \mc{E})$, which is equivalent to finding the stability number of $\mc{G}$'s complementary graph $\bar{\mc{G}} = (\mc{N}, \bar{\mc{E}})$. Let $\omega$ be the maximum clique number of graph $\mc{G}$, then we can formulate the maximum clique problem as the following MIP, which finds the stability number on graph $\bar{\mc{G}}$:
\begin{subequations}\label{eq: max_clique_mip}
\begin{alignat}{4}
\omega = \max ~~&\sum_{i=1}^n x_i\\
\st ~~&x_i+x_j\leq1 &~~&\forall (i,j)\in\bar{\mc{E}}\\
&x_i\in\{0,1\} &~~&\forall i = 1,...,n.
\end{alignat}
\end{subequations}

Let $A$ be the adjacency matrix of $\mc{G}$, then we have $A = Q - \bar{A}$, where $Q = \mathbf{e} \mathbf{e}^\top - I$, $\bar{A}$ is the adjacency matrix of $\bar{\mc{G}}$. Applying this relationship to the COP model in Corollary 2.4 of \cite{de2002approximation}, we obtain the following COP model to calculate the maximum clique number of $\mc{G}$:
\begin{subequations}\label{eq: max_clique_cop}
\begin{alignat}{4}
\omega = \min ~~&\lambda\\
\st~~&\lambda (\mathbf{e} \mathbf{e}^\top - A) - \mathbf{e} \mathbf{e}^\top = Y\label{eq: max_clique_cop_1}\\
&Y\in\mc{C}_{n}.
\end{alignat}
\end{subequations}

In our experiment we use 10 max-clique problem instances from the second DIMACS challenge. We compare the following ways of solving the COP \eqref{eq: max_clique_cop}:

(1) Approximately solve the COP with SDP, as shown in Section \ref{ec: sdp_approx}. This is the method suggested by \cite{de2002approximation}. We use Mosek 9.1.2 to solve those SDP approximations.

(2) Exactly solve the COP with the cutting plane algorithm of Section \ref{ch: cut_plane}. 

Notice that we can strengthen the master problem in our cutting plane algorithm by providing bounds for $Y$ in the initialization stage. Since the maximum clique number $\omega$ cannot exceed the number of total nodes $|\mc{N}|$,  and elements of $\mathbf{e} \mathbf{e}^\top$ are all 1's while elements of $A$ are either 0 or 1, from constraints \eqref{eq: max_clique_cop_1} we have that the elements of $Y$ should be in the range of $[-1, |\mc{N}| - 1]$.

We present the results of the algorithmic comparison in Table \ref{tab: algo_max_clique}, where we list the number of nodes $|\mc{N}|$, number of edges $|\mc{E}|$, and the maximum clique number of the graph $\omega$ for each instance. For the computational performance of solving the SDP approximation via Mosek, we list the objective (``Obj"), optimality gap (``Gap", compared with the true $\omega$), and the computational time (``Time"). For the performance of our cutting plane algorithm we list the computational time and the number of iterations needed for convergence. There is no need to list the objectives because the cutting plane method always converges to an exact solution.
{\small
\begin{table}[h]
\centering
\caption{Algorithm Comparison for Maximum Clique COP Model}
\label{tab: algo_max_clique}
\begin{tabular}{lrrrrrrrrrr}
\hline
              & \multicolumn{1}{l}{}       & \multicolumn{1}{l}{}       & \multicolumn{1}{l}{}      &  & \multicolumn{3}{c}{Mosek}                                                                        & \multicolumn{1}{c}{} & \multicolumn{2}{c}{Cutting plane}                                \\
Instance      & \multicolumn{1}{c}{$|\mc{N}|$} & \multicolumn{1}{c}{$|\mc{E}|$} & \multicolumn{1}{c}{$\omega$} &  & \multicolumn{1}{c}{Obj} & \multicolumn{1}{c}{Gap(\%)} & \multicolumn{1}{c}{Time(sec)} &                      & \multicolumn{1}{c}{Time(sec)}         & \multicolumn{1}{c}{\#Iter}        \\ \cline{1-4} \cline{6-8} \cline{10-11} 
c-fat200-1    & 200                        & 1534                       & 12                        &  & 12                            & 0                                & 566.81                        &                      & 13.87                                 & 2                                       \\
c-fat200-2    & 200                        & 3235                       & 24                        &  & 24                            & 0                                & 638.72                        &                      & 18.90                                 & 2                                       \\
c-fat200-5    & 200                        & 8473                       & 58                        &  & 60.35                         & 3.89                            & 606.33                        &                      & 12.19                                 & 2                                       \\
hamming6-2$^a$     & 64                         & 1824                       & 32                        &  & 32                            & 0                                & 1.51                    &      &6.05                      & 2 \\
hamming6-4$^a$    & 64                         & 704                        & 4                         &  & 4                             & 0                                & 1.59                   &       & 1.55                     &4      \\
johnson8-2-4  & 28                         & 210                        & 4                         &  & 4                             & 0                                & 0.20                          &                      & 9.53                                  & 2                                       \\
johnson8-4-4  & 70                         & 1855                       & 14                        &  & 14                            & 0                                & 2.47                          &                      & 11.82                                 & 2                                       \\
johnson16-2-4$^{a,b}$ & 120                        & 5460                       & 8                         &  & 8                             & 0                                & 31.88                 &        &      62.75                &   2                            \\
keller4       & 171                        & 9435                       & 11                        &  & 13.47                         & 18.34                           & 426.16                        &     &-                 &-                            \\
MANN\_a9      & 45                         & 918                        & 16                        &  & 17.48                         & 8.47                            & 0.45                          &                      & 547.62                                & 2                                       \\ \hline
\multicolumn{9}{l}{\footnotesize $^a$ Obtained by setting $q = \bar{\omega}$ in the separation problem, see text for explanation. }\\
\multicolumn{9}{l}{\footnotesize $^b$ Obtained by early termination of the separation problem. }
\end{tabular}
\end{table}}

When solving instances ``hamming6-2", ``hamming6-4" and ``johnson16-2-4" with cutting planes, we encounter some very hard separation problems that take very long time to solve. To speedup the process, we use instead a strengthened version of the separation problem with $q = \bar{\omega}$ in constraint \eqref{eq: separation_c} \citep{anstreicher2020testing}, where $\bar{\omega}$ is the current master problem solution for $\omega$. In addition, even after our enhancement of $q = \bar{\omega}$, the instance ``johnson16-2-4" still has a hard separation problem which achieves a nonzero lower bound early on (thus proves that the matrix is not copositive), but cannot converge after an extended period of time. In this case we set a time limit of 1 minute to help the separation problem stop early. Finally, when solving the instance ``keller4" we have a very hard separation after a few iterations. CPLEX fails to find a feasible solution for this separation problem after an extended period of time, and we had to stop the solution process because of high memory usage. However, we can still get useful information from the master problem objective, as it always provides an upper bound for the COP objective. When we stopped the algorithm for ``keller4" the master problem objective was already 11, which equals the correct value of $\omega$ is better than the bound provided by the SDP approximation (13.47).

From the results we can observe that for some instances, the SDP approximation fails to provide the correct maximum clique number. Also, in certain instances such as ``c-fat200-1", ``c-fat200-2" and ``c-fat200-3", the cutting plane algorithm is faster than the SDP approximation.

We also compare our algorithm with the simplicial partition method of \cite{bundfuss2009adaptive}, which we believe is the only exact algorithm for general linear COPs in the literature. \cite{bundfuss2009adaptive} also solve the maximum clique instances from the second DIMICS challenge. They report that their computation time for ``johnson8-2-4" and ``hamming6-4" are respectively 1 minutes 33 seconds and 57 minutes 52 seconds. For all the other instances, their algorithm produces only poor bounds within two hours. Therefore, the performance of our algorithm is better than theirs in all test instances. 

Notice that the cutting plane algorithm terminates in very few iterations for almost all the test instances. It is not generally the case with the cutting plane algorithm when solving other COP problems.  One reason for the small number of iterations could be the use of a strong formulation for the maximum clique problem. For example, if we use the weaker COP formulation \eqref{eq: complex_max_clique_cop} below, then the cutting plane algorithm takes longer to terminate: the simplest instance (in terms of the number of nodes and edges) ``johnson8-2-4" now takes 200.64 seconds and 690 iterations:
\begin{subequations}\label{eq: complex_max_clique_cop}
\begin{alignat}{4}
\omega = \min ~~&\lambda\\
\st~~&\lambda I + \sum_{(i,j)\in\bar{\mc{E}}}x_{ij} E_{ij}  - \mathbf{e} \mathbf{e}^\top = Y\\
&Y\in\mc{C}_{n}
\end{alignat}
\end{subequations}
where $E_{ij}\in\mb{R}^{n\times n}$ is a matrix with ones at $i$th row and $j$th column and  $j$th row and $i$th column, and with zeros for all other positions. 
\section{Plots for Trends of Revenue and Profit}\label{ec: trend_graph}
\begin{figure*}[h]
\centering
\begin{subfigure}[t]{0.35\textwidth}
\centering
\includegraphics[width = \textwidth]{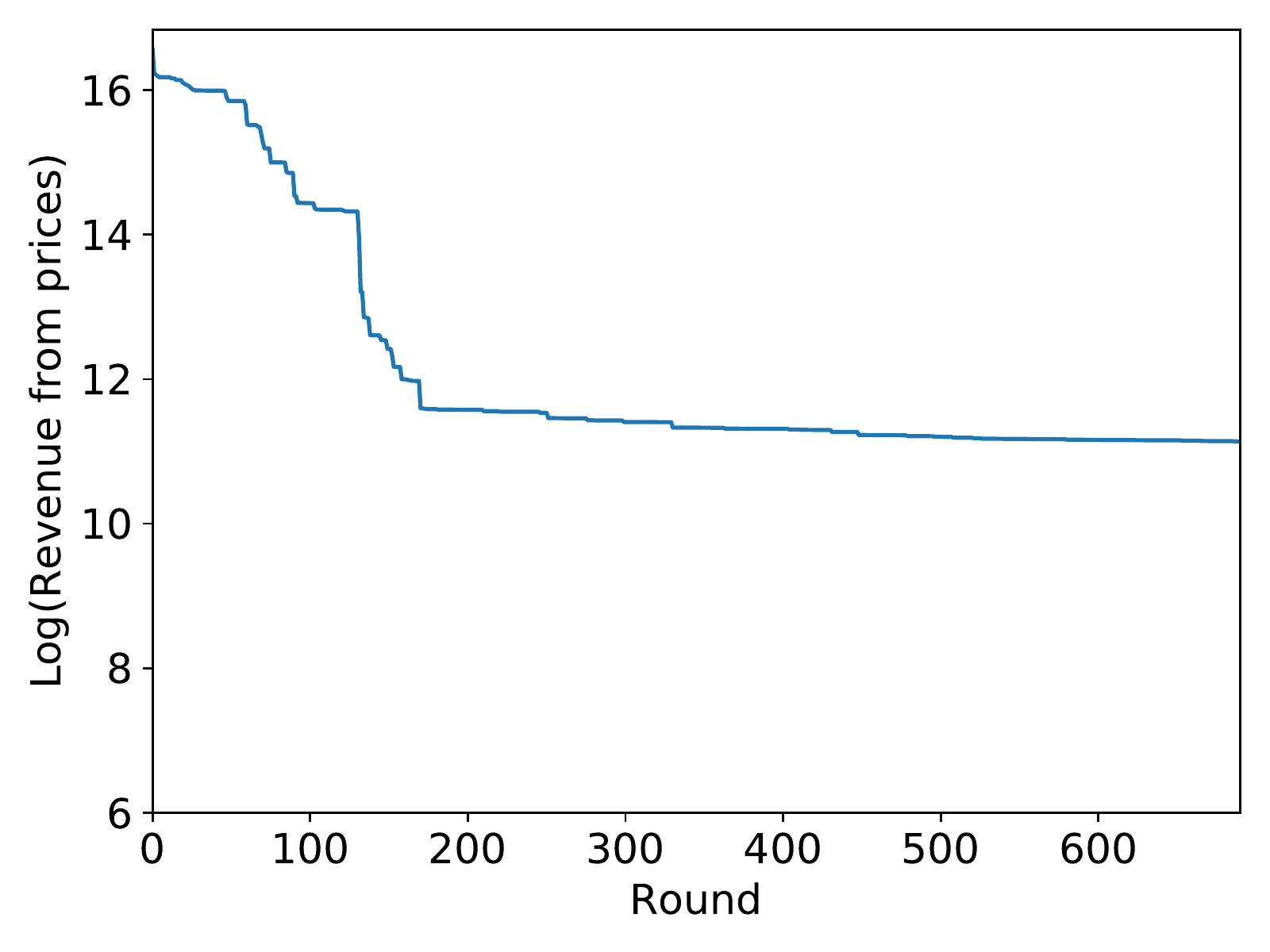}
\caption{}
\end{subfigure}
~
\begin{subfigure}[t]{0.35\textwidth}
\centering
\includegraphics[width = \textwidth]{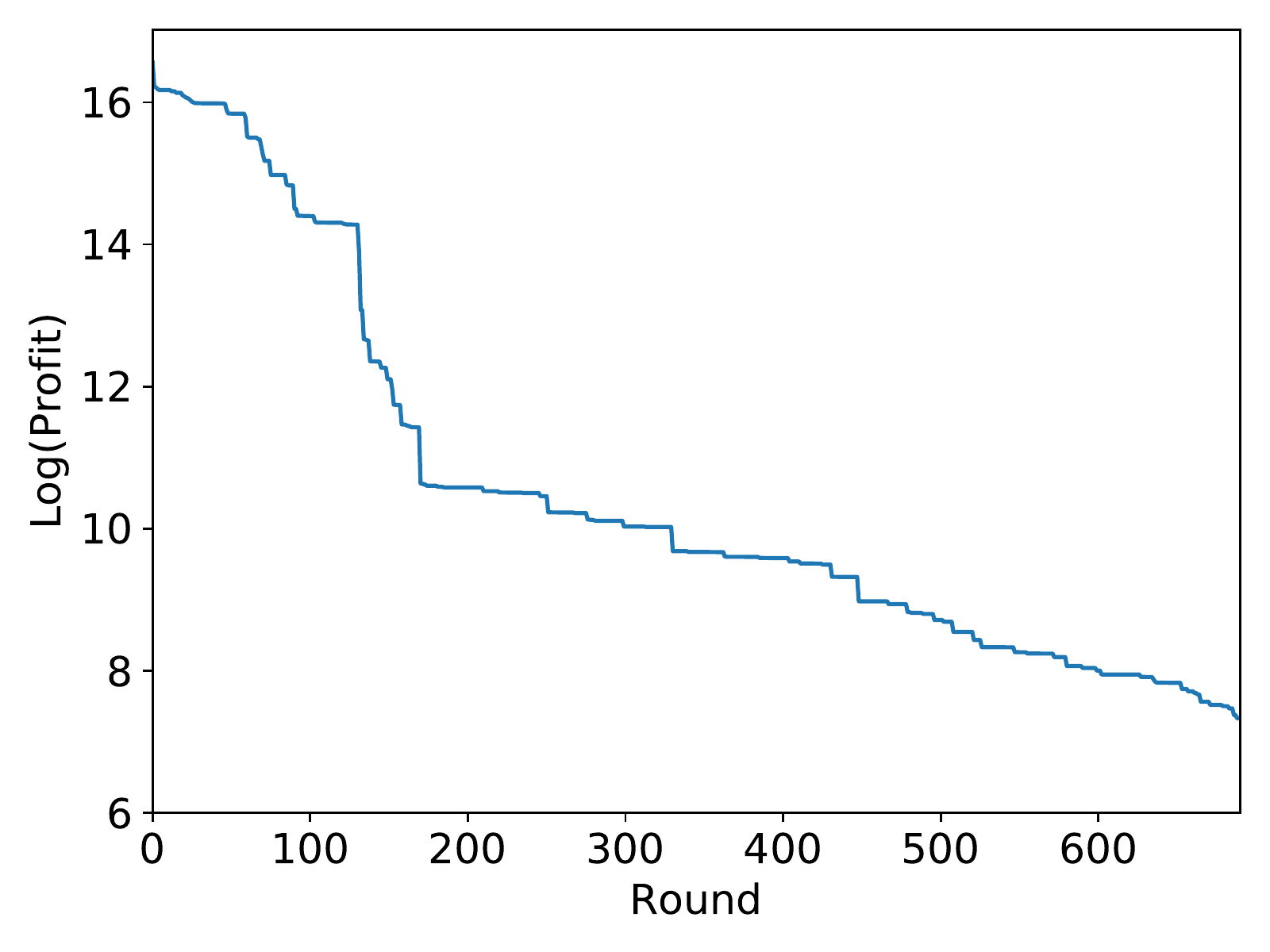}
\caption{}
\end{subfigure}
~
\begin{subfigure}[t]{0.35\textwidth}
\centering
\includegraphics[width = \textwidth]{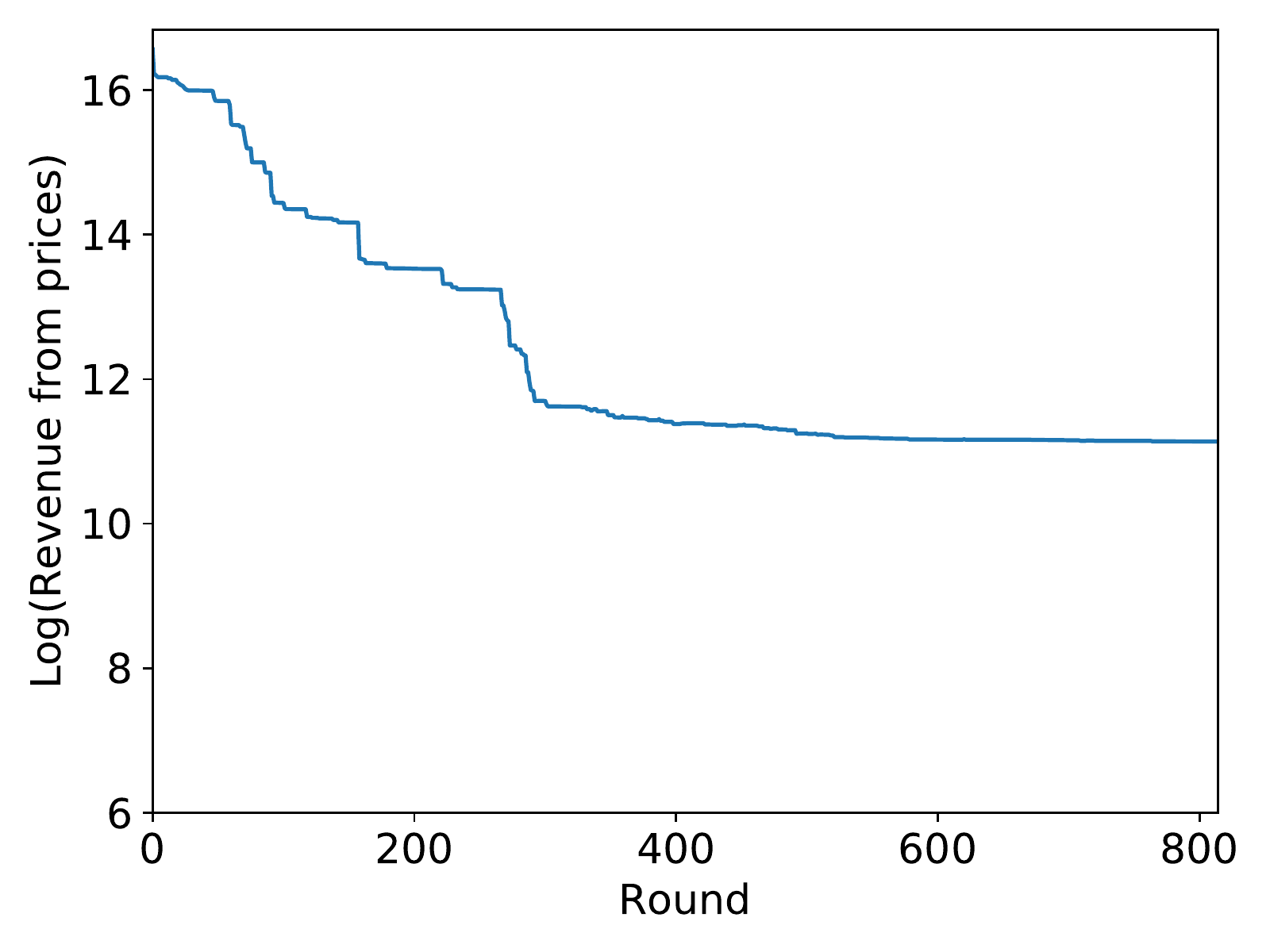}
\caption{}
\end{subfigure}
~
\begin{subfigure}[t]{0.35\textwidth}
\centering
\includegraphics[width = \textwidth]{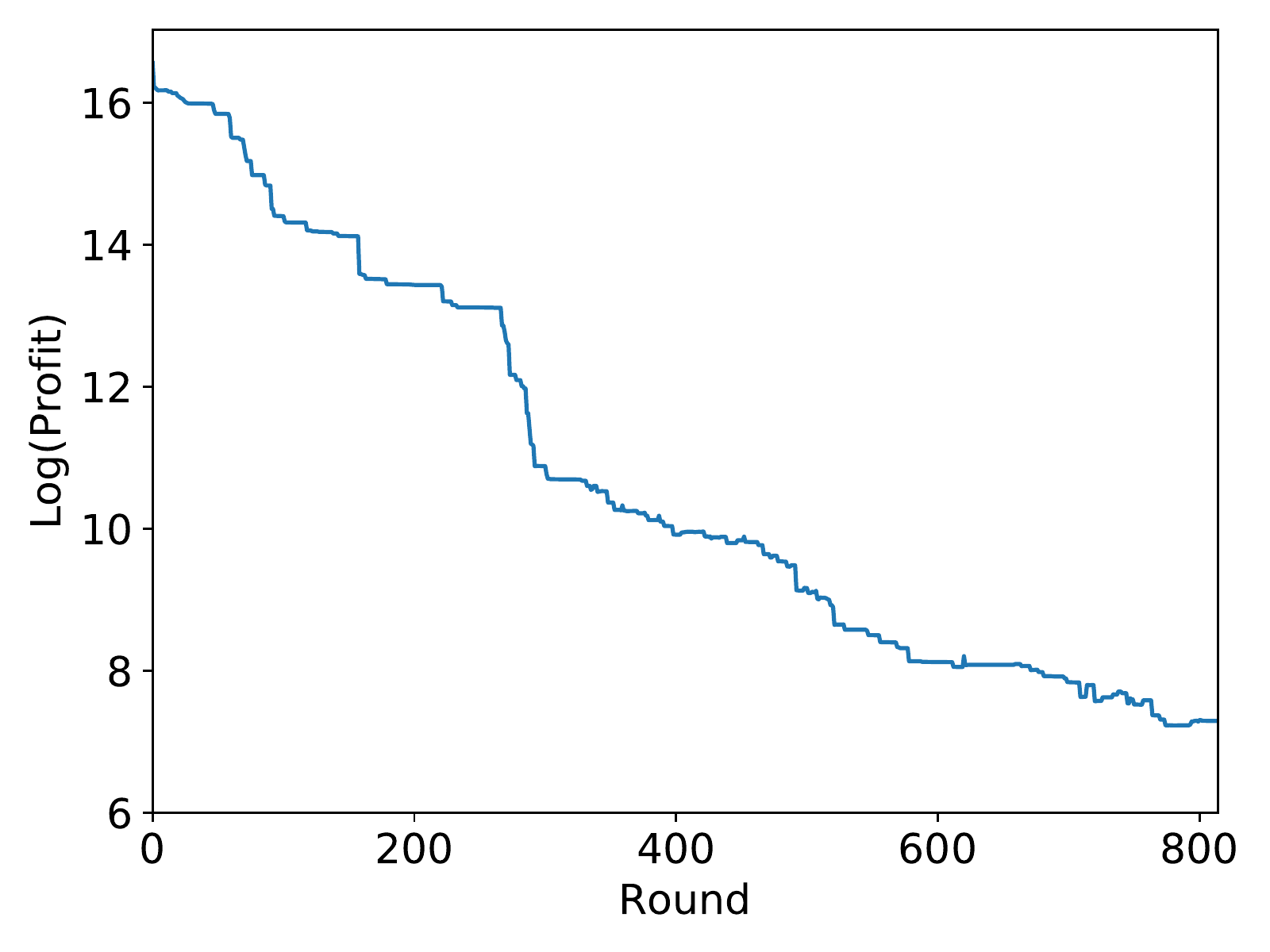}
\caption{}
\end{subfigure}
\caption{\small Trends of the cutting plane algorithm for the first UC instance in (a) revenue of CDP, (b) profit of CDP, (c) revenue of RCDP, (d) profit of RCDP.}
\label{fig: trends}
\end{figure*}

\begin{figure*}[h]
\centering
\begin{subfigure}[t]{0.35\textwidth}
\centering
\includegraphics[width = \textwidth]{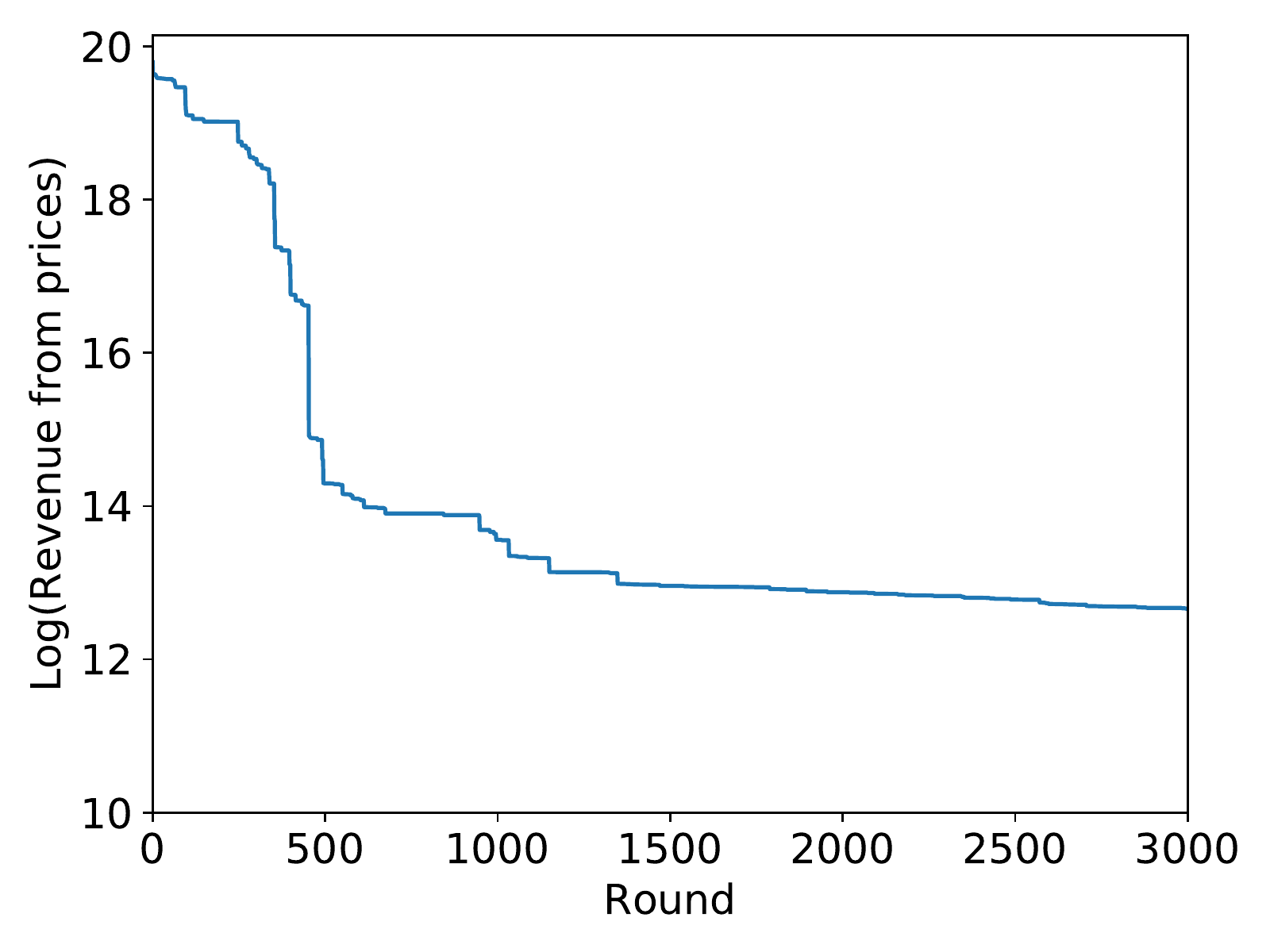}
\caption{}
\end{subfigure}
~
\begin{subfigure}[t]{0.35\textwidth}
\centering
\includegraphics[width = \textwidth]{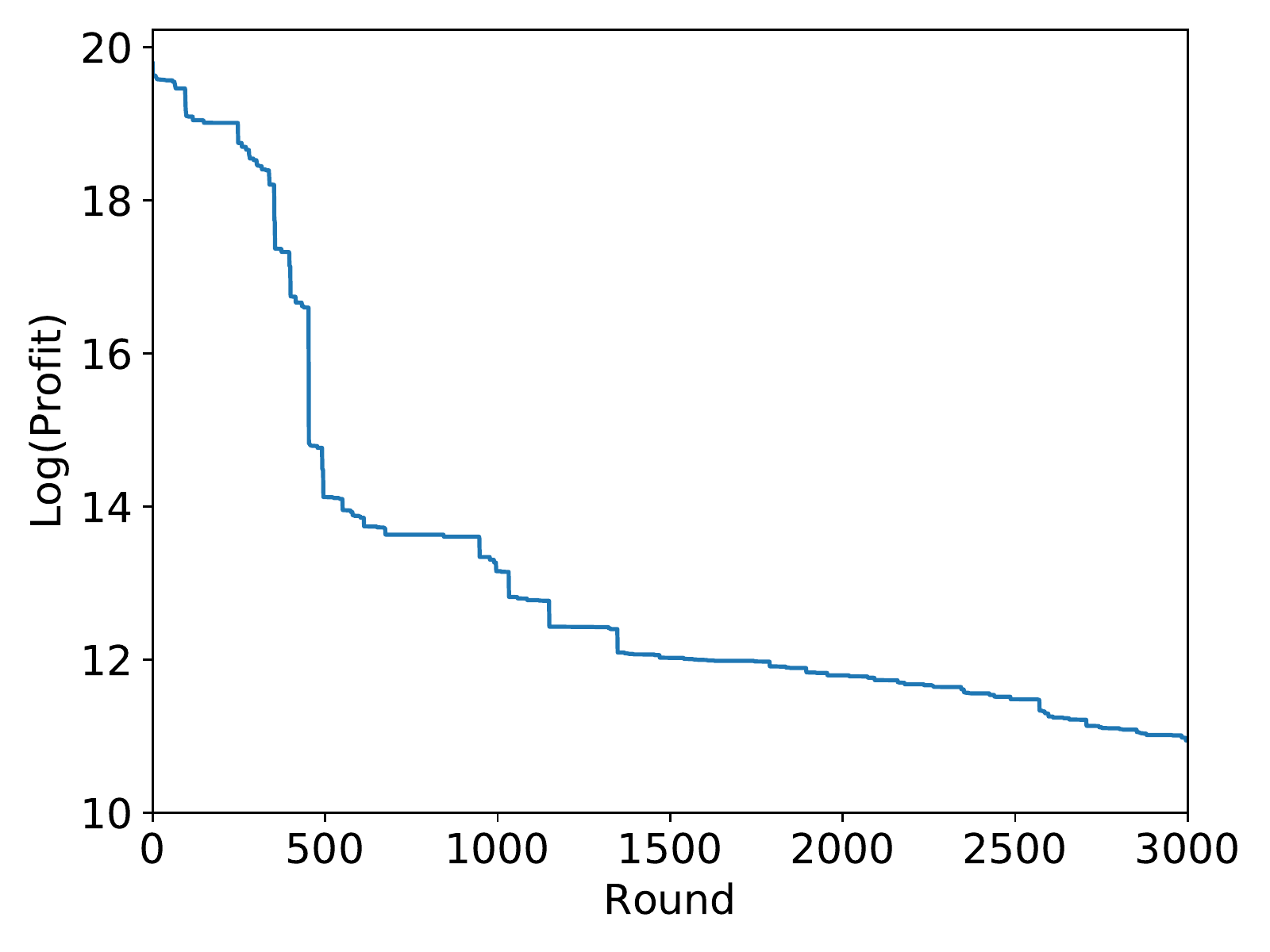}
\caption{}
\end{subfigure}
~
\begin{subfigure}[t]{0.35\textwidth}
\centering
\includegraphics[width = \textwidth]{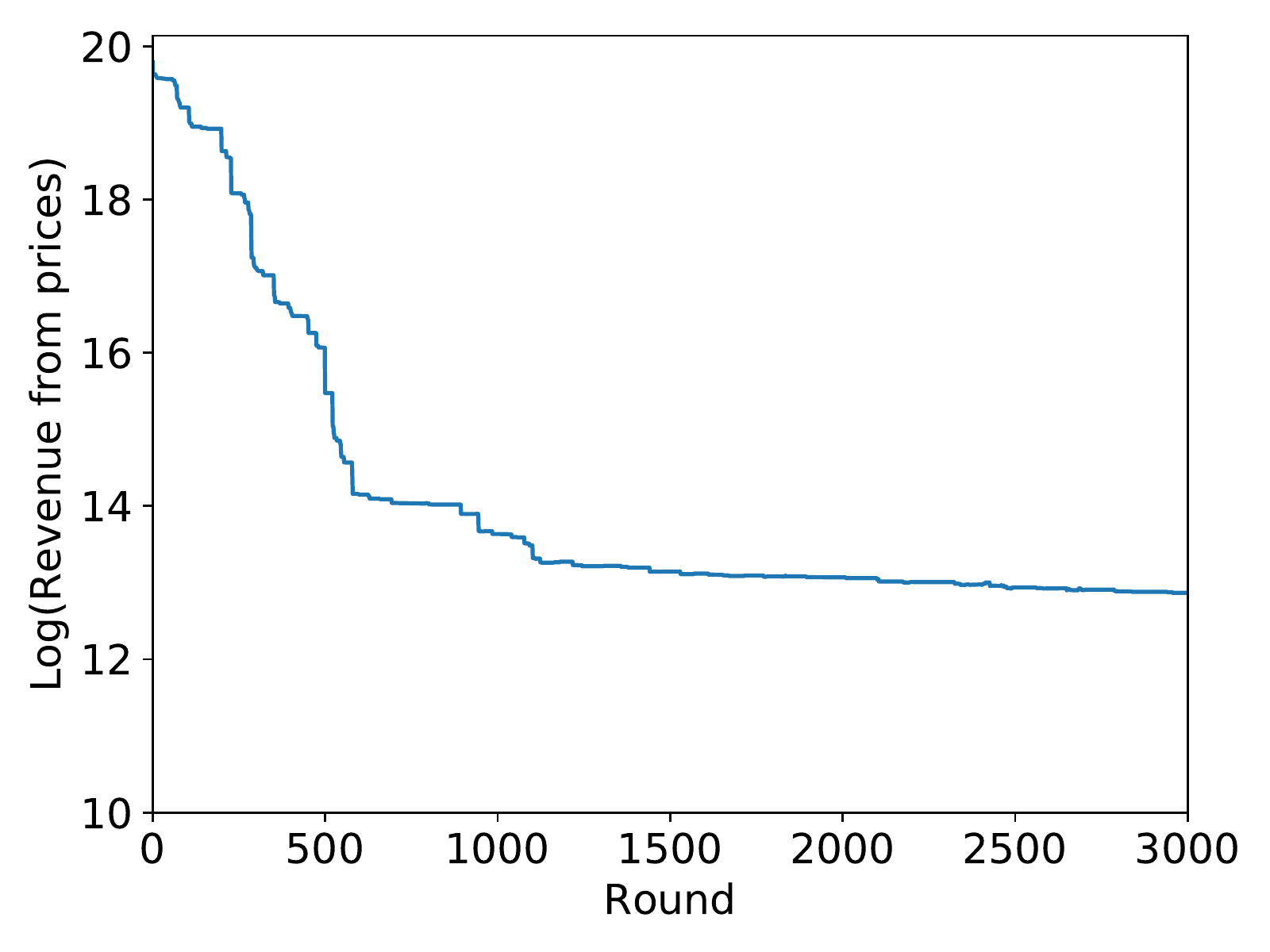}
\caption{}
\end{subfigure}
~
\begin{subfigure}[t]{0.35\textwidth}
\centering
\includegraphics[width = \textwidth]{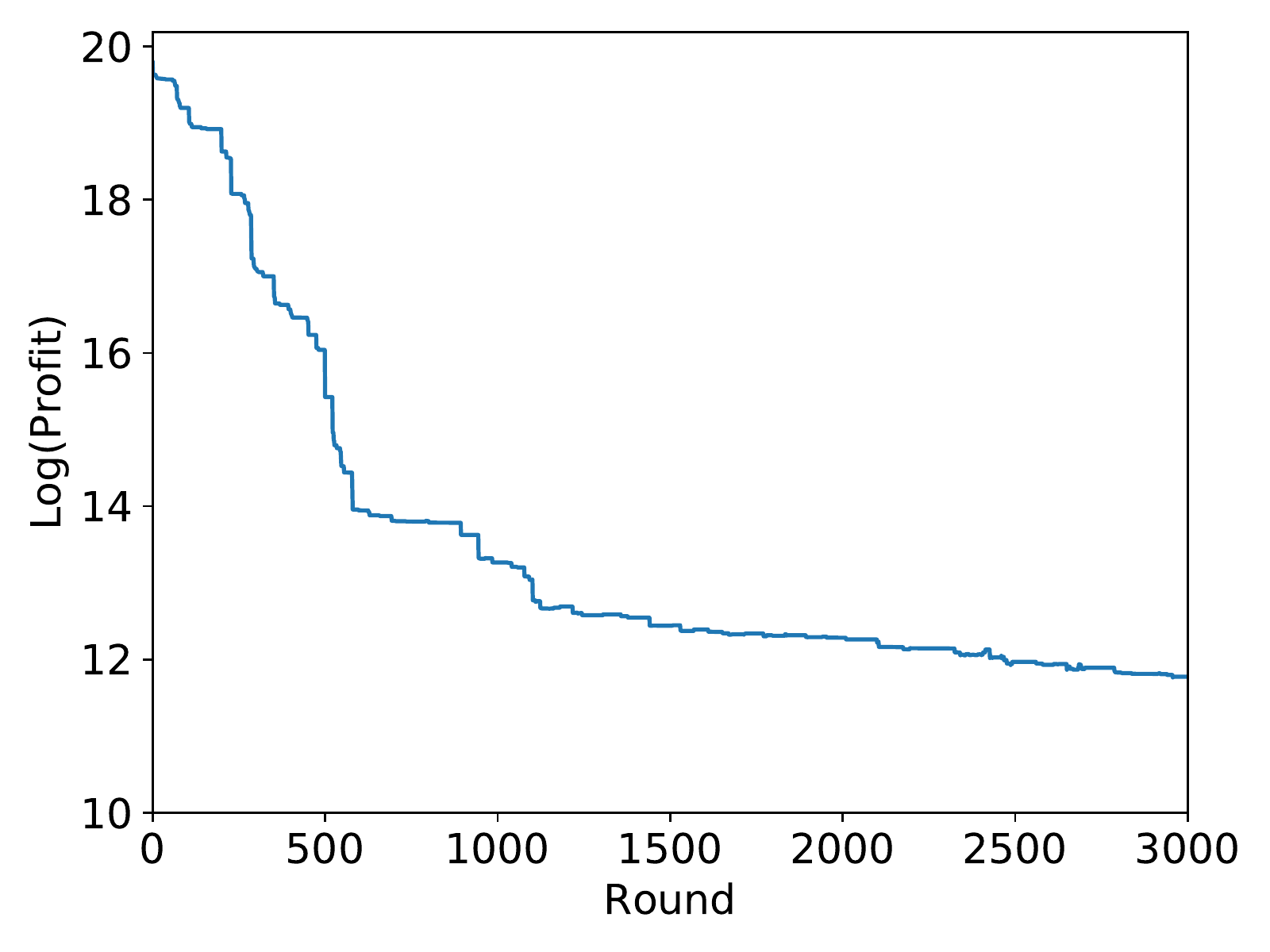}
\caption{}
\end{subfigure}
\caption{\small Trends of the cutting plane algorithm for the second UC instance in (a) revenue of CDP, (b) profit of CDP, (c) revenue of RCDP, (d) profit of RCDP.}
\label{fig: trends2}
\end{figure*}

In Figure \ref{fig: trends} we show the revenue and profit trends in the cutting plane algorithm for the first UC instance. We observe that for both CDP and RCDP, the revenue and profit are monotonically decreasing. Similar trends are observed for the second UC instance in Figure \ref{fig: trends2} (except for very small increases for a few data points). Those trends show a tendency of convergence, which means that the algorithm still produces meaningful results even if we terminate it early.
